%% file: potential_July2026.tex
\documentclass[11pt]{amsart}

\usepackage[hmargin=0.8in,height=8.8in]{geometry}
\usepackage{amssymb,amsthm, times}
\usepackage{delarray,verbatim}
\usepackage{srcltx}
\usepackage{ifpdf}
\ifpdf
\usepackage[pdftex]{graphicx}
 \else
\usepackage[dvips]{graphicx}
 \fi

\usepackage{caption}
\usepackage{xcolor}
\usepackage{mathrsfs}
\definecolor{webgreen}{rgb}{0,.5,0}
\definecolor{webbrown}{rgb}{.8,0,0}
\definecolor{emphcolor}{rgb}{0.95,0.95,0.95}

\usepackage{hyperref}
\hypersetup{%
          colorlinks=true,
          linkcolor=webbrown,
          filecolor=webbrown,
          citecolor=webgreen,
          breaklinks=true}
\ifpdf \hypersetup{
            pdfstartview=FitH, 
            bookmarksopen=true,
            bookmarksnumbered=true
} \else \hypersetup{dvips} \fi
\usepackage{enumitem}
\usepackage{comment}
\usepackage[normalem]{ulem}
\usepackage[colorinlistoftodos,prependcaption,textsize=small]{todonotes}
\setuptodonotes{color=blue!20} 
\usepackage{tikz}

\numberwithin{equation}{section}

\newtheorem{proposition}{Proposition}[section]
\newtheorem{corollary}{Corollary}[section]
\newtheorem{remark}{Remark}[section]
\newtheorem{lemma}{Lemma}[section]

\newtheorem{assump}{Assumption}
\newtheorem{definition}{Definition}[section]
\numberwithin{remark}{section} \numberwithin{proposition}{section}
\numberwithin{corollary}{section}

\newcommand {\R}{\mathbb{R}}

\newcommand {\F}{\mathcal{F}}

\newcommand {\p}{\mathbb{P}}

\newcommand {\E}{\mathbb{E}}
\newcommand {\LL}{\mathcal{L}}

\newcommand{\eps}{\varepsilon}
\newcommand{\diff}{{\rm d}}

\newcommand{\conn}{\quad\text{and}\quad}
\newcommand{\1}{\mbox{1}\hspace{-0.25em}\mbox{l}}
\newcommand{\lev}{L\'{e}vy }

\newcommand{\goes}{\rightarrow \infty}
\newcommand{\eq}{\mathbf{e_q}}

\title{A potential-theoretic approach to optimal stopping in a spectrally  L\'{e}vy Model}
\author[M. Egami]{Masahiko Egami}
\address[M. Egami]{Graduate School of Economics,
Kyoto University, Sakyo-Ku, Kyoto, 606-8501, Japan}
\email{egami.masahiko.8x@kyoto-u.ac.jp}
\thanks{
}
\author[T. Koike]{Tomohiro Koike}
\address[T. Koike]{Graduate School of Economics,
Kyoto University, Sakyo-Ku, Kyoto, 606-8501, Japan}
\email{koike.tomohiro.4e@kyoto-u.ac.jp}

\date{}
\begin{document}
\begin{abstract}
We establish a systematic solution method for optimal stopping problems of spectrally negative \lev processes. 
Our approach relies essentially on potential theory, that is, on the analysis of superharmonic and subharmonic functions and their analytic properties, including the maximum principle.
Using these mathematical results,  we not only derive necessary and sufficient conditions of optimality for a broad class of reward functions, but also develop a method to tackle general problems in a direct and constructive way (without pre-specifying the solution form).   
To reinforce the latter point, we also present several examples with complex solution structures that illustrate the effectiveness of our approach, including continuation regions with multiple connected components.
\end{abstract}

\maketitle \noindent \small{\textbf{Key words:} Optimal stopping; spectrally negative \lev processes; potential theory; 
maximum principle;
Gerber-Shiu function. }\\

\section{Introduction}\label{sec:introduction}
 \subsection{The problem}
This paper investigates the optimal stopping problem for a spectrally negative \lev process, a class of real-valued \lev processes with no positive jumps.
Let the spectrally negative \lev process $X=\{X_t;t\geq 0\}$ represent the state variable defined on the probability space $(\Omega, \F, \p)$, where $\Omega$ is the set of all possible outcomes, and $\p$ is a probability measure defined on $\F$. 
 For $x \in \mathbb{R}$, let $\p_x$ denote the law of $X$ started at $x$.
We write $\E_x$ for the associated expectation.
We denote by $\mathbb{F}=\{\F_t\}_{t\ge 0}$ the filtration with respect to which $X$ is adapted, assuming that the usual conditions hold.

For a spectrally negative \lev process, the Laplace exponent $\psi$  is given by
\begin{align}\label{eq:psi}
\psi(\theta)=-\gamma\theta+\frac{1}{2}\sigma^2\theta^2+\int_{(-\infty,0)}(e^{\theta x}-1-\theta x\1_{(-1,0)}(x))\Pi(\diff x),
\end{align}
where $\gamma \in \R$, $\sigma \geq 0$, and $\Pi$ is a measure concentrated on $(-\infty, 0)$ satisfying
$\int_{(-\infty,0)}(1\wedge x^2)\Pi(\diff x)<\infty$. It is well-known that $\psi$ is zero at the origin
and has a right-continuous inverse: for $q \geq 0$,
\[\Phi(q) :=\sup\{\lambda \geq 0: \psi(\lambda)=q\}.\]
Moreover, it is known that  $X_t \goes $ as $t \goes$ almost surely if and only if $\psi'(0+)>0$, oscillates if and only if $\psi'(0+)=0$ and $X_t \to -\infty $ as $t \goes$ almost surely if $\psi'(0+)<0$.
The jumps of the process have a finite mean $\int_{(-\infty, 0)} |x| \Pi(\diff x)< \infty$ and there is no diffusion component $\sigma = 0$ if and only if the paths have bounded variation. Then we may rewrite \eqref{eq:psi} as
\begin{align*}
\psi(\theta) = \delta \theta + \int_{(-\infty, 0)} \Big( e^{\theta x}-1\Big) \Pi(\diff x),\quad \text{where}
\end{align*}
\begin{equation}\label{eq:delta}
\delta := - \left(\gamma + \int_{(-1, 0)} x \Pi(\diff x)\right)
\end{equation}
 is the drift coefficient.

Let $\LL$ be the infinitesimal generator of $X$, which  is given by
\begin{align}\label{eq:operator-L}
\LL f(x) = -\gamma f'(x) +\frac{1}{2}\sigma^2 f''(x)+ \int_{(-\infty, 0)}\left[ f(x+y)-f(x)-y \1_{(-1,0)}(y)f'(x)\right]\Pi(\diff y).
\end{align}

Let $q \ge 0$ be a constant and $g$ be a non-negative Borel function.
We denote by
\begin{align}
    v(x) = \sup_{\tau \in \mathcal{T}} \E_x [e^{-q\tau} g(X_\tau)], \ x \in \R \label{eq:value}
\end{align}
the value function of the optimal stopping problem with reward function $g$ and discount rate $q$, where the supremum is taken over the class $\mathcal{T}$ of all $\mathbb{F}$-stopping times.
We define the stopping region $\Gamma$ and the continuation region $C$ as follows:
\begin{align}\label{eq:Gamma-C}
     \Gamma= \{x \in \R: v(x)=g(x) \} \conn
   C = \{x \in \R: v(x)>g(x) \}.
\end{align}

We make the following assumptions. 
\begin{assump}\label{assumption-jumpsize}
\begin{enumerate}[label=(\roman*)]
     \item  $X$ has paths of bounded variation with $\delta>0$, or paths of unbounded variation with $\sigma>0$.
     
     \item  If $X$ has paths of bounded variation, then
$g$ is upper semi-continuous and right continuous with left limits (RCLL), and belongs to
$C^1(\R\setminus F)$, where $F$ is a finite set.
Furthermore, for every $x\in F$, the one-sided derivatives
$g'(x\pm)$
exist and are finite.
\item  If $X$ has paths of unbounded variation, then
$g$ is upper semi-continuous and continuous, and belongs to
 $C^2(\R\setminus F)$, where $F$ is a finite set.
Furthermore, for every $x\in F$, the one-sided derivatives
 $g'(x\pm)$ and $g''(x\pm)$
exist and are finite.
     \item For any $\eps>0$ and $x\in\R$, there exists $\delta>0$ such that
\[
\int_{(\eps,\infty)}
\sup_{u\in[x-\delta,x+\delta]}
|g(u-z)|\,\Pi(-\diff z)
<\infty.
\] 
   . 
     \item The \lev measure $\Pi$ satisfies $\int_{(-\infty, -1)} |x| \Pi(\diff x)<\infty.$

\item $\Pi$ has no atoms.

\end{enumerate}
\end{assump}
To simplify the notation, we simply say that $X$ is of bounded variation (resp.\ unbounded variation) when it is of bounded variation with $\delta>0$ (resp.\ of unbounded variation with $\sigma>0$), unless confusion may arise.
Under (ii) or (iii), $g$ is locally bounded, since it is either RCLL or continuous.
Note that, under (ii) and (iii), $g(X)$ is RCLL. In the unbounded-variation case, this follows immediately from the continuity of $g$. In the bounded-variation case, since $X$ has no positive jumps, upward crossings occur continuously, and the right-continuity of $g$ ensures the right-continuity of $g(X)$.
Assumptions (ii), (iii), (iv), and (v) are imposed to carry out the analysis via the generator $\LL$.
Under these conditions, $\LL g$ is finite and continuous on $\R \setminus F$.
These properties can be established by applying the dominated convergence theorem.
Assumption (v) is equivalent to the process $X_t$ having a finite mean for each fixed $t$ and is imposed also in \cite{avram2015gerber}, \cite{avram2020w}, and \cite{biffis2010note}.
Assumption (vi) is adopted to ensure the smoothness of the scale function, which is introduced in Section \ref{sec:scale_functions}.

Assumption (iv) is in fact equivalent to the following condition, although the latter may at first glance appear to be weaker;
\begin{enumerate}
    \item[(iv')] For each $x\in\R$, there exist $\eps>0$ and $\delta>0$ such that
\[
\int_{(\eps,\infty)}
\sup_{u\in[x-\delta,x+\delta]}
|g(u-z)|\,\Pi(-\diff z)
<\infty.
\] 
\end{enumerate}
The implication from (iv') to (iv) follows because $\Pi$ is finite on sets bounded away from zero and $g$ is locally bounded.

\subsection{Literature review and research motivation}\label{sub:literature}
Optimal stopping problems for \lev processes have been extensively studied. One strand of research specifies the reward function in explicit form and derives constructive solutions for the associated stopping problems. For instance, \cite{mordecki}, \cite{alili-kyp}, and \cite{christensen2009note} investigated the McKean optimal stopping problem (see Section \ref{subsection:onesided}), while \cite{kyprianou2005novikov} and \cite{novikov2007solution} focused on problems with power-type reward functions. The papers of \cite{egami-yamazaki2013} and \cite{Egami-Yamazaki-2011} consider an optimal alarm problem in the context of capital adequacy management, where the reward function is decreasing and negative in the negative region and identically zero in the positive region. 
In these works, the expected reward is computed explicitly, and the set maximizing it is identified.
As a consequence, the reward function must be chosen so that the expected payoff admits an explicit representation.

A more recent strand of research has shifted attention toward studying optimal stopping problems with broader classes of reward functions, without specifying their functional form in advance. One line of research makes use of 
the \textit{averaging problem}, a term introduced in \cite{surya2007approach}. The pioneering works in this direction include \cite{surya2007approach}, \cite{deligiannidis2009optimal}, and \cite{mordecki2007optimal}. Specifically, \cite{surya2007approach} and \cite{deligiannidis2009optimal} focus on the one-sided case and characterize the optimal boundary as the root of the averaging function, which is the solution of the averaging problem. Lemma 3.1 in \cite{surya2007approach} provides a general and constructive method for computing this function. \cite{mordecki2007optimal} characterizes the solution to optimal stopping problems for a broad class of Hunt processes using the representation theory of excessive functions. By leveraging the Green kernel for Lévy processes expressed in terms of the process maxima, this approach recovers results analogous to those obtained in \cite{surya2007approach} and \cite{deligiannidis2009optimal}. Furthermore, \cite{christensen2013optimal} characterized the solution of optimal stopping problems for general Hunt processes by the averaging problem. This result is an extension of  \cite{surya2007approach}, \cite{deligiannidis2009optimal} and \cite{mordecki2007optimal}.  They also provide a method for constructing the solution of the averaging problem for \lev processes, linear diffusions and continuous-time Markov chains. 
It enables one to solve the problem constructively \textit{in the one-sided setting}. \cite{mordecki2016optimal} refines the averaging problem approach by reducing the domain required for the averaging function, thus obtaining explicit solutions for a broader class of reward functions. \cite{mordecki2021two} presents a two-sided verification  with the averaging problem and provides an example in which a two-sided solution arises.

However, two notable features common to the averaging problem approach  are the following: (1) the form of the solution must be specified in advance, and (2) it seems difficult to systematically address problems more complex than the one-sided case. To the best of our knowledge, when the problem is not one-sided, no general constructive method for solving the averaging problem is known. Indeed, cases where the solution can be identified as one-sided are rather limited, while in practice the more interesting problems are often two-sided or beyond. Moreover, although the aforementioned approaches have the advantage of being applicable to general \lev processes, it should be noted that  explicit solutions are typically not available because the value function obtained by the authors' method is represented as a function of the running minimum of the process, which is often not known explicitly.

In comparison with the literature, here is our contribution: motivated by this open question, we establish a systematic approach that does \textit{not require the prior specification of the solution form}, thereby avoiding the need for a guess-and-verify procedure. While our focus on spectrally negative \lev processes may entail some loss of generality, we present a powerful and constructive method applicable to a broad class of reward functions, covering one-sided, two-sided, and multiple mixed cases. This is particularly helpful in cases such as Sections~\ref{subsection:twosided}, \ref{subsection:discontinuous}, and \ref{subsec:bimdal},
where the conventional guess-and-verify method becomes intractable.
Moreover, it is worth noting that within our framework, fluctuation identities are not the primary analytical tool; instead, they play a secondary role when deriving explicit solutions in specific examples.

A key element enabling our approach is effective use of potential-theoretic techniques.
We heavily make use of  analytical properties of sub- or superharmonic functions to construct the solution in Section \ref{sec:sufficient}.
The relative positions of the reward function $g$ and the candidate value function differ depending on whether $g$ is sub- or superharmonic (Lemmas~ \ref{lemma:sub-dominance}, \ref{lemma:majorant} and \ref{lemma:monotonicity}).
These results are based on the maximum principle introduced in Proposition \ref{prop:subsuperMP}.
Using these results, we construct the true value function, which dominates $g$ and  is superharmonic on the entire set.

It should be noted that the application of potential theory to optimal stopping problems has been explored in several studies. 
One of the earliest contributions in this direction is \cite{salminen1985}, which characterized the solution to the optimal stopping problem for one-dimensional diffusion using the Martin representation of the value function. As discussed above, \cite{mordecki2007optimal} and \cite{christensen2013optimal} are  significant works in this line for \lev processes. More recently, this methodology has been extended to multi-dimensional diffusions. For example, \cite{christensen2018multidimensional} and \cite{christensen2019optimal} derived integral equations examining the stopping region through the Martin boundary theory.

In particular, in Section \ref{sec:sufficient}, \textit{the smooth Gerber–Shiu function} introduced in \cite{avram2015gerber} plays a crucial role. 
This function is a harmonic function on an interval $(a, \infty)$ that is smooth at the boundary $a$. 
The candidate value function can be expressed with the smooth Gerber–Shiu function; see \eqref{eq:expectedreward} and \eqref{eq:expectedreward2}.
For an overview of Gerber–Shiu theory, see Chapter 10 of \cite{Kyprianou_2014} and \cite{kyprianou2013gerber}.

Finally, we note another line of research focused on the log-concavity of the reward function. The representative works are \cite{hsiau2014logconcave} and  \cite{lin2019one}. They show that the log-concavity and monotonicity of the reward function imply the existence of a one-sided solution under
general random walks in discrete time and Lévy processes in continuous time.  
The relationships between these works and our paper are described in Corollary \ref{coro:logconcave}.

The rest of the paper is organized as follows. In Section \ref{sec:preliminary}, we provide the mathematical preliminaries necessary for this paper. We mainly review the basic properties of spectrally negative \lev processes, the fundamental notions of sub- or superharmonic functions, including the maximum principle, and the smooth Gerber–Shiu function. 
Section \ref{sec:martin} constitutes one of the main contributions of this paper, as it identifies 
where given candidate of the value function is sub- or superharmonic.
The main tool in the proof  is Meyer -- It\^o formula.  
We use this representation to derive the necessary condition in Section \ref{sec:necessary} and to check the superharmonicity of the candidate value function in Section \ref{sec:sufficient}.
Section \ref{sec:necessary} demonstrates that the continuous fit condition (resp. the smooth fit condition) at the boundary points of  the continuation region is necessary when $X$ has paths of bounded variation (resp. unbounded variation).
Section \ref{sec:sufficient} gives a sufficient condition for optimality and how to construct the solution, which is the main result of this paper.
Section \ref{sec:example} presents various examples in which we solve specific problems following the procedure presented in Section \ref{sec:sufficient}.

\section{Mathematical tools}\label{sec:preliminary}
\subsection{Spectrally negative  \lev processes and their scale functions} \label{sec:scale_functions}
For every spectrally negative \lev  process, there exists a $q$-scale function
\(
W^{(q)}: \R \rightarrow \R\) for \( q\ge 0,
\)
that is continuous, strictly increasing on $[0,\infty)$  and $0$ on $(-\infty,0)$. It is
uniquely determined by
\begin{eqnarray*}
\int_0^\infty e^{-\beta x} W^{(q)}(x) \diff x = \frac{1}{\psi(\beta)-q}; \quad \beta > \Phi(q).
\end{eqnarray*}
When \( q = 0 \), we simply write \( W \) instead of \( W^{(0)} \).
Under  Assumption \ref{assumption-jumpsize}-(vi), $W^{(q)}$ belongs to $C^1(0,\infty)$, and to $C^2(0,\infty)$ if $X$ has paths of unbounded variation and $\sigma > 0$; see  \cite{Chan_2009}.
For $q > 0$, or for $q = 0$ with $\psi'(0+) < 0$, the scale function increases exponentially:
\begin{align}
   \lim_{x \goes} W^{(q)} (x) /e^{\Phi(q) x} = 1/\psi'(\Phi(q)). \label{eq:limit-W}
\end{align}
For $q=0$ and $\psi'(0+)>0$, $ \lim_{x \goes} W (x)  = 1/\psi'(0+)$.
It follows that (see \cite{Hubalek_Kyprianou_2009})
\begin{align}
    W^{(q)}(0+)&=  \label{eq:W0}
    \begin{cases}
        1/ \psi'(0+) &\text{if} \ X \ \text{has paths of bounded variation},  \\
        0  &\text{if} \ X \ \text{has paths of unbounded variation} ,
    \end{cases}\\
    W^{(q)\prime}(0+)&= \label{eq:Wprime0}
    \begin{cases}
        (q+\Pi(-\infty, 0))/(\psi'(0+))^2 &\text{if} \ X \ \text{has paths of bounded variation},  \\
        2/\sigma^2  &\text{if} \ X \ \text{has paths of unbounded variation},
    \end{cases}
\end{align}
where $W^{(q)\prime}(0+)$ may be infinity.
For a comprehensive account of the scale function, we refer the reader to
\cite{Bertoin_1996,Bertoin_1997, Kyprianou_2014, Kyprianou_Surya_2007}.  
See also  \cite{Egami_Yamazaki_2010_2, Surya_2008} about methods for numerically computing the
scale function.

\subsection{Superharmonic, subharmonic, and harmonic Functions}\label{subsection:potential}
Let $\{H_t \}_{t \geq 0} $ be the Markov kernel of $X$; for $x \in \mathbb{R}$ and a Borel set $A$ and for a Borel function $f$,
\[
H_t(x,A)=\p_x (X_t \in A) \conn H_t f (x) = \E_x [f(X_t)],
\]
respectively.  Moreover, we extend this notation from deterministic time $t$ to stopping times $\tau$:
\(H_\tau(x,A)=\p_x (X_\tau \in A).\)
For any Borel set $A$, we define the hitting time
\begin{align*}
    T_A := \inf \{t \ge 0: X_t \in A\} \quad \text{ and write} \quad H_A:=H_{T_A}.
\end{align*}
Moreover, for simplicity, we write
\( T_r:= T_{(r, \infty)}=\inf\{t\ge 0: X_t > r\}\) and 
\(T_\ell^{-}:=T_{(-\infty, \ell)}=\inf\{t\ge 0: X_t <\ell\}\)
for $\ell$ and $r$ in $\R$.
In this vein, we will use, throughout this paper, the following notation for hitting times of $(-\infty, a)$ and $(-\infty, a) \cup (b, \infty)$ for $a, b \in \mathbb{R}$:
\begin{align*}
    &H_a := H_{(-\infty, a)}, \quad H_{a,b} := H_{(-\infty, a) \cup (b, \infty)}.
\end{align*}
For $a = -\infty$, we naturally set $H_{-\infty, b} := H_{(b, \infty)}$.
For $q \geq 0$, $H^q_t$ is defined to indicate
\( H^q_t (x,A) = \p_x (X_t \in A; t<\eq),\)
where $\eq$ is a random variable independent of $X$ and follows an exponential distribution with rate $q$.
Similarly, $H_\tau f$, $H^q_t f$, 
$H^q_\tau f$ are defined in the same manner.
In particular, for $x \in \mathbb{R}$, $H^q_\tau f$ is defined by
\begin{align*}
H^q_\tau f(x) = \E_x \bigl[e^{-q \tau} f(X_\tau)\mathbf{1}_{\{\tau<\infty\}}\bigr],
\end{align*}
and will be used repeatedly later. 


A non-negative measurable function \( u \) is said to be \textit{\( q \)-superharmonic} (resp. \textit{\( q \)-subharmonic}) on an open set $G$ if each open subset $A \subset G$ whose closure is compact in $G$,
\(   u \geq  H_{A^\mathrm{c}}^q u \quad (\text{resp. } \ u \leq  H_{A^\mathrm{c}}^q u ). \)
A \textit{q-harmonic function} is defined as a function that is both $q$-superharmonic and $q$-subharmonic.
A function is said to be \textit{strictly superharmonic} if it is superharmonic but not harmonic. Similarly, we define a \textit{strictly subharmonic} function as one that is subharmonic but not harmonic.

\subsection{The smooth Gerber–Shiu function}\label{sec:G-S}
For \( a \in \mathbb{R} \), let \( h_a^{(q)} \) be a $q$-harmonic function on \( (a, \infty) \) satisfying the boundary condition \( h_a^{(q)}(x+) = g(x-) \) for \( x \leq a \)
(and \( h_a^{(q)\prime}(a+) = g'(a-) \) when \( X \) has paths of unbounded variation). This function is called the \textit{smooth Gerber–Shiu function}, and it is known that it can be represented in terms of \( W \)  as follows (\cite{avram2015gerber}): for every $x \geq a$, 
\begin{align}
  h_a^{(q)}(x) = g(a-) + g'(a -)(x-a) - \int_0^{x-a} W^{(q)}(x-a-y) J_a(y) \, \diff y, \label{eq:GerberShiu} 
\end{align}
where \( J_a \) is given by, for each $x>0$, 
\begin{align*}
    J_a(x) =& g'(a-) \psi'(0+) - q(g'(a-)x + g(a-)) \\
           &+ \int_{(x,\infty)} \left( g(x+a-z) - g(a-) + g'(a-)(z-x) \right) \Pi(-\diff z).
\end{align*}

When \( q = 0 \), we simply write \( h_a \) instead of \( h_a^{(0)} \).
The same convention applies to $h_a^{(q)\prime}$.
Under Assumption \ref{assumption-jumpsize}, the function $J_a(x)$ is continuous. Consequently, $h_a^{(q)}$ is of class $C^1$ on $(a,\infty)$ if $W^{(q)}$ is of class $C^1$, and $ h_a^{(q)}$ is of class $C^2$ on $(a,\infty)$ if $W^{(q)}$ is of class $C^2$ by an argument similar to that in the proof of Lemma 5.7 in \cite{avram2015gerber}.
Moreover, the mapping $a \mapsto  h_a^{(q)}(x)$ is left-continuous, which is proved in Proposition \ref{prop:h_a-left-cont} in Appendix \ref{sec:h_a-left-cont}.


Note that $h_a$ satisfies Assumption \ref{assumption-jumpsize} (iv).
Since $h_a = g$ on $(-\infty, a)$, we have, for sufficiently large $\eps>0$,
\[ \int_{(\eps,\infty)}
\sup_{u\in[x-\delta,x+\delta]}
|h_a(u-z)|\,\Pi(-\diff z) = \int_{(\eps,\infty)}
\sup_{u\in[x-\delta,x+\delta]}
|g(u-z)|\,\Pi(-\diff z)\]
The right-hand side is finite by Assumption \ref{assumption-jumpsize} (iv); and hence the left-hand side is also finite. Therefore, $h_a$ satisfies Assumption \ref{assumption-jumpsize} (iv'), which is equivalent to (iv).

The asymptotic behavior at infinity is given by (Lemma 5.7 in \cite{avram2015gerber}):
\begin{align}
    \lim_{x \to \infty}& \frac{h_a^{(q)}(x)}{W^{(q)}(x-a)} = \kappa(a), \quad \text{where} \label{eq:kappa} \\
    \kappa(a) = &\frac{\sigma^2}{2}g'(a-) + \frac{q}{\Phi(q)}g(a-) 
    - \int_0^\infty \diff x\, e^{-\Phi(q)x} \label{eq:kappa-def}\\
    &\int_{(x, \infty)} \left(g(x+a-z) - g(a)\right) \Pi(-\diff z), \nonumber
\end{align}
with the convention that the second term on the right-hand side  of \eqref{eq:kappa-def} is interpreted as
$\psi'(0+)=\lim_{q\to0} q/\Phi(q)$ when $q=0$ and $X$ drifts to infinity.
We refer to $\kappa$ as the $\kappa$-function associated with $g$.
It follows from Fubini's theorem and a direct calculation that
\begin{align}\label{eq:derivative-kappa}
e^{-\Phi(q)a}\kappa(b)-e^{-\Phi(q)a}\kappa(a) =  \int_a^b e^{-\Phi(q)x}(q-\LL g) (x) \diff x.
\end{align}

The $\kappa$-function provides a convenient way to identify the boundary of the continuation region in one-sided optimal stopping problems (see Corollaries \ref{coro:one-sided} and \ref{coro:logconcave} ).

One obtains (Propositions 5.4 and 5.5 in \cite{avram2015gerber})
\begin{align}
   H^{(q)}_a g(x) &= h^{(q)}_a(x) - W^{(q)}(x-a) \kappa(a), \label{eq:expectedreward}\\
  H^{(q)}_{a,b} g(x) &= h^{(q)}_a(x)+W^{(q)}(x-a)\frac{g(b)-h^{(q)}_a(b)}{W^{(q)}(b-a)}.\label{eq:expectedreward2}
\end{align}
In view of the smoothness properties of $h^{(q)}_a$ and $\kappa$ established above,
both $H^{(q)}_a g$ and $H^{(q)}_{a,b} g$ are of class $C^1$ on $(a,\infty)$ and $(a,b)$, respectively, 
if $W^{(q)}$ is of class $C^1$.
Moreover, they are of class $C^2$ on the same domains if $W^{(q)}$ is of class $C^2$.

\subsection{The maximum principle}\label{subsection:mp}
We prove a maximum principle for subharmonic and superharmonic functions. This principle is heavily used in Section \ref{sec:sufficient}.

\begin{proposition}\label{prop:subsuperMP}
Assume that $f$ is superharmonic (resp. subharmonic) on an interval $(a,b)$.
Then there is no $x \in (a,b)$ such that $f(x) < f(y)$ (resp. $f(x)>f(y)$) for every $y \in (-\infty, a]\cup \{b\}$. 
\end{proposition}
\begin{proof}
    We only prove the superharmonic case, since the subharmonic case follows analogously.
    Assume to the contrary that there exists  $x \in (a,b)$ such that $f(x) < f(y)$  for all $y \in (-\infty, a]\cup \{b\}$. 
    Since the support of the exit distribution
    $\mathbb P^x(X_{T_{(a,b)^{\mathrm c}}}\in \mathrm dy)$
    is contained in
    $(-\infty,a]\cup\{b\}$
    under the spectrally negative assumption, we have
    \[ H_{a,b}f(x) = \int f(y) \p_x [X_{T_{(a,b)^{\mathrm{c}}}} \in \diff y] >f(x),\]
    which contradicts the assumption that
     $f$ is superharmonic on $(a,b)$.
\end{proof}

\section{A representing measure of expected reward functions}\label{sec:martin}
Our approach, presented in Section \ref{sec:sufficient}, requires identifying whether the candidate value function is superharmonic or subharmonic. The purpose of this section is therefore to establish a method for determining the regions where a given function $p$ is $q$-superharmonic or $q$-subharmonic. Throughout the subsequent analysis, $p$ will be taken to be an expected reward function associated with a stopping strategy, namely,
$p=H_{A} g$,
which represents the expected reward obtained by stopping upon the first entry into a given set $A$.
The key ingredient in the proof of Proposition \ref{prop:Riesz} is the Meyer--It\^o formula.
We use the following notations;
\begin{align*}
    \Delta f(x) := f(x+)-f(x-),\conn
    \Delta X_t := X_t -X_{t-}.
\end{align*} 

\begin{proposition}\label{prop:Riesz}
    Assume that $X$ has paths of bounded variation  (resp. unbounded variation). 
    Assume that a non-negative function $p$ is of class $C^1$ (resp.\ $C^2$) except at a finite set $F^p$, that
\(
\Delta p^{(k)}(a)
\)
is finite for $k=0,1,2$ and every $a\in F^p$, and that $p$ is RCLL (resp.\ continuous).
Let $\mu$ be a signed measure such that 
\[
\mu(dx)=
\begin{cases}
(q-\LL)p(x)\,dx, & x\notin F^p,\\[1ex]
-\delta  \Delta p(x) \ \left(\text{resp.} \ -\dfrac{\sigma^{2}}{2}
\Delta p'(x), \right)
& x\in F^p,
\end{cases}
\]

Then, \(\mu|_{S}\ge0\)
on an open set $S$ if and only if
$p$ is $q$-superharmonic on $S$.
\end{proposition}

\begin{remark}
In the later sections, we apply Proposition~\ref{prop:Riesz} with
$p=g$, $p=H_ag$, and $p=H_{a,b}g$.
We briefly verify that these functions satisfy its assumptions.
For $p=g$, this follows immediately from Assumption~\ref{assumption-jumpsize}~(ii) and (iii).
For $p=H_ag$ and $p=H_{a,b}g$, the functions are of class
$C^1$ (resp.~$C^2$) on $(a,\infty)$ and $(a,b)$, respectively, when $X$ has paths of bounded variation (resp.~unbounded variation). Moreover, they are RCLL (resp.~continuous) at the boundary points $a$ and $b$. Hence, they also satisfy the assumptions of Proposition~\ref{prop:Riesz}.
\end{remark}


\begin{lemma} \label{lemma:meyer-ito}
Let $p$ be as defined in Proposition~\ref{prop:Riesz}.
    \begin{enumerate}
        \item Assume that $X$ has paths of unbounded variation. Then, we have
        \begin{align}\label{eq:unbounded-meyer-ito}
            e^{-qt}&p(X_t) - p(X_0) = \int_0^t e^{-qs}(\LL-q)p(X_{s-})\diff s   \\
             &+\frac{1}{2}\sum_{a \in F^p} \Delta p'(a) \int_0^t e^{-qs} \diff L_s^a + M_t, \nonumber
        \end{align}
        where $L_t^a$ is the local time at $a$ of $X$ and $(M_t)_{t \geq 0}$ is some local martingale. 
        \item Assume that $X$ has paths of bounded variation. Then, we have\begin{align}\label{eq:bounded-meyer-ito}
            e^{-qt}&p(X_t) - p(X_0) = \int_0^t e^{-qs}(\LL-q)p(X_{s-})\diff s \\
            &+\sum_{a \in F^p} \Delta p(a) \int_0^t e^{-qs} \diff N_s^{a\uparrow} + M_t, \nonumber
        \end{align}
        where $N_t^{a\uparrow}: = \# \{s\in (0,t]: X_{s-}=X_s=a \} $ and $(M_t)_{t \geq 0}$ is some local martingale.
    \end{enumerate}
\end{lemma}
\begin{proof}
    We only show the case where $q=0$. The case where $q>0$ can be proved by using the result for $q=0$ and the following integration by parts (Corollary 2  of Theorem II. 20 in \cite{protter2012stochastic}):
    \[ 
    e^{-qt}p(X_t)-p(X_0)
= \int_0^t e^{-qs} \diff p(X_s)
-
q\int_0^t e^{-qs}p(X_{s-}) \diff s.\]

First we prove (1).
Assume  $X$ has paths of unbounded variation.
    Since $p$ is a difference of two convex functions (see Problem 3.6.20 in \cite{karatzas}),
    the Meyer--It\^o formula (Theorem IV. 70 in \cite{protter2012stochastic}) leads to 
    \begin{align}\label{eq:meyer-ito}
        p(X_t)-p(X_0) = &\int_0^t p'(X_{s-})\diff X_{s} + \sum_{0 \leq s \leq t}[p(X_s)-p(X_{s-})-p'(X_{s-})\Delta X_s]   \\
        &+\frac{1}{2}\int_{-\infty}^\infty D^2 p (\diff a)L_t^a,\nonumber
    \end{align}
    where $D^2$ is the second derivative in the sense of distribution  and it is given by
    \[D^2 p = p'' + \sum_{a \in F^p}[p'(a+)-p'(a-)]. \]
    By the occupation density formula (Corollary 1 of Theorem IV. 70 in \cite{protter2012stochastic}), we compute the last term in \eqref{eq:meyer-ito} as 
    \begin{align}\label{eq:occupation}
        \frac{1}{2}\int_{-\infty}^\infty D^2 p (\diff a)L_t^a
        &= \frac{1}{2}\int_{-\infty}^\infty p''(a) L_t^a \diff a + \frac{1}{2}\sum_{a \in F^p} \Delta p'(a) L_t^a,\\
        &= \frac{\sigma^2}{2}\int_0^t p''(X_{s-})\diff s + \frac{1}{2}\sum_{a \in F^p} \Delta p'(a) L_t^a.\nonumber
    \end{align}
    
    By It\^o -\lev decomposition of $X$ (Theorem 2.1 in \cite{Kyprianou_2014}), we compute the first term in \eqref{eq:meyer-ito} as 
    \begin{align}\label{eq:ito-lev}
        \int_{0}^t p'(X_{s-})\diff X_{s} =& -\gamma \int_0^t p'(X_{s-})\diff s + \sum_{0\leq s \leq t}p'(X_{s-})\Delta X_s \1_{\{\Delta X_s \leq -1\}} \\ \nonumber
        &+ (\text{a local martingale})
    \end{align}
Combining \eqref{eq:occupation} and \eqref{eq:ito-lev} with \eqref{eq:meyer-ito}, we obtain
\begin{align} \label{eq:meyer-ito2}
    p(X_t)-p(X_0) 
    &= \int_0^t  \left(-\gamma p'(X_{s-})+\frac{\sigma^2}{2}p''(X_{s-}) \right)\\
    &+ \underbrace{\sum_{0 \leq s \leq t}\left[ p(X_s)-p(X_{s-})-\1_{\{\Delta X_s>-1\}} p'(X_{s-})\Delta X_s \right]}_{(\bigstar)}\nonumber \\
    &+ \frac{1}{2} \sum_{a \in F^p} \Delta p'(a) L_t^a + (\text{a local martingale}). \nonumber
\end{align}
    Finally, we compute ($\bigstar$). Let $N(\diff s, \diff x)$ be the Poisson random measure of $X$. 
    \begin{align}\label{eq:compensation}
        (\bigstar) 
        =& \int_0^t  \int_ {-\infty}^0 \Big( p(X_{s-}+x)-p(X_{s-})-\1_{\{x>-1\}}x p'(X_{s-}) \Big)  N(\diff s, \diff x),\\
        =& \int_0^t \int_{-\infty}^0 \Big( p(X_{s-}+x)-p(X_{s-})-\1_{\{x>-1\}}xp'(x) \Big) \Pi(\diff x) \diff s \nonumber \\
        &+(\text{a martingale}), \nonumber
    \end{align}
    where the second equality follows from the compensation formula (Theorem 4.4 in \cite{Kyprianou_2014}). Substituting this into \eqref{eq:meyer-ito2} and using \eqref{eq:operator-L}, we obtain \eqref{eq:unbounded-meyer-ito}.

    We prove (2).
    Assume that $X$ has paths of bounded variation. Note that $p$ may have discontinuities, and hence we cannot apply the Meyer--It\^o formula directly.
    We use the following notation: 
    \begin{align}        
        p(X)_{s-} &:= \lim_{s' \uparrow s} p(X_{s'}) . \label{eq:p-limit}
    \end{align}
    Note that
\(p(X)_{s-}=p(X_{s-})\)
unless $X$ creeps continuously through a discontinuity point of $p$ at time $s$. 
    We show that the following identity holds for any $p$ satisfying the following conditions: $p\in C^1(\mathbb{R}\setminus F)$, $\Delta p^{(k)}(a)$ is finite for $k=0,1,2$ and $a\in F$, and $p$ is RCLL. 
    \begin{align}\label{eq:meyer-ito-bdd}
        p(X_t)-p(X_0) = \int_0^t p'(X_{s-})\diff X_s + \sum_{0 \leq s \leq t} (p(X_s)-p(X)_{s-}-p'(X_{s-})\Delta X_s).
    \end{align}
    Note that $p$ can be decomposed as
\[
p=f+\sum_{a\in F}\Delta p(a)\1_{[a, \infty)},
\]
where $f$ is the difference of two convex functions. Since the Meyer--It\^o formula applies to $f$, \eqref{eq:meyer-ito-bdd} holds for $f$. Therefore, by linearity, it suffices to prove \eqref{eq:meyer-ito-bdd} for $\1_{[a, \infty)}$, which that \eqref{eq:meyer-ito-bdd} holds for any $p$.
Indeed, we have, for any fixed $a$,
\begin{align*}
  & \1_{[a, \infty)}(X_t)-\1_{[a, \infty)}(X_0)\\
   &=
    \begin{cases}
        0, \quad \text{if} \ (X_t<a \text{ and } X_0<a)
\quad\text{or}\quad
(X_t\ge a \text{ and } X_0\ge a), \\
        1, \quad \text{otherwise},
    \end{cases}\\
    &= \sum_{0 \leq s \leq t} \left( \1_{[a, \infty)}(X_s)-(\1_{[a, \infty)}(X))_{s-} \right),
\end{align*}
 and hence, \eqref{eq:meyer-ito-bdd} holds for $\1_{[a,\infty)}$.
With the notation in \eqref{eq:p-limit} with $p= \1_{[a,\infty)}$, note that the term $ \1_{[a,\infty)}(X_s)-\big(\1_{[a,\infty)}(X)\big)_{s-} $
takes the value \(1\) when the path crosses \(a\) from left to right at time \(s\), which occurs through continuous creeping, and the value \(-1\) when the path crosses \(a\) from right to left, in which case the crossing occurs by a jump.

Using It\^o--\lev decomposition for $X$ as in the unbounded variation case, it follows from \eqref{eq:meyer-ito-bdd} that
\begin{align}\label{eq:meyer-ito2-bdd}
     p(X_t)&-p(X_0) \\
     = &\int_0^t -\gamma p'(X_{s-})\diff s + \underbrace{\sum_{0 \leq s \leq t} (p(X_s)-p(X)_{s-}-p'(X_{s-})\Delta X_s\1_{(-1 < \Delta X_s < 0)})}_{(\bigstar \bigstar)} \nonumber \\
     &+(\text{a local martingale}). \nonumber
\end{align}
Finally, we compute $(\bigstar \bigstar)$.
Noting that \(p(X)_{s-}=p(X_{s-})\)
holds at every jump time $s$, i.e., whenever $\Delta X_s<0$, we have
\begin{align*}
    (\bigstar \bigstar) 
   = \sum_{0 \leq s \leq t}& (p(X_s)-p(X)_{s-}-p'(X_{s-})\Delta X_s\1_{(-1 < \Delta X_s < 0)})\1_{(\Delta X_s<0)} \\
   &+ \sum_{0 \leq s \leq t} (p(X_s)-p(X)_{s-}) \1_{(\Delta X_s = 0)}.
\end{align*} 
By using the compensation formula as in \eqref{eq:compensation}, the first term of $(\bigstar \bigstar)$ can be expressed, up to a local martingale, as
\[\int_0^t \int_{(-\infty, 0)} p(X_{s-}+x)-p(X_{s-})-p'(X_{s-})x\1_{(-1<x<0)}  \Pi(\diff x) \diff s.\]
Recalling the definition of $N_t^{a\uparrow}$ in Lemma \ref{lemma:meyer-ito} (2), the second term of $(\bigstar \bigstar)$ is given by
\[   \sum_{a \in F^p} \Delta p(a) N_t^{a \uparrow}.\]
 Substituting this into \eqref{eq:meyer-ito2-bdd} and using \eqref{eq:delta} and \eqref{eq:operator-L}, we obtain \eqref{eq:bounded-meyer-ito}.
\end{proof}

\begin{proof}[Proof of Proposition \ref{prop:Riesz}]

We only present the proof for the unbounded variation case, since the proof for the bounded variation case is analogous.
Assume that $\mu|_S \geq 0$ for an open $S$. We show that $p$ is
 $q$-superharmonic on $S$.
Take any $A \subset S$ with compact closure. It follows from Lemma \ref{lemma:meyer-ito} 
\begin{align}\label{eq:doob-mayer}
e^{-q{(t\wedge\sigma_A})}p(X_{t\wedge\sigma_A})-p(X_0)
=
A_{t\wedge\sigma_A}
+
M_{t\wedge\sigma_A},
\end{align}
where $\sigma_A$ is the first exit time from $A$, and
\[ A_t: = \int_0^t e^{-qs}(\LL-q)p(X_{s-})\diff s + \frac{1}{2}\sum_{a \in F^p} \Delta p'(a) \int_0^t e^{-qs} \diff L_s^a.\]
The proof is complete if we show that
$e^{-q{(t\wedge\sigma_A})}p(X_{t\wedge\sigma_A})$ is a supermartingale.
Since $p$ is non-negative and $A_{t\wedge\sigma_A}\le 0$ by the assumption
$\mu|_S\ge0$ ,
$M_{t\wedge\sigma_A}$ is a local martingale bounded below,
and hence it is a supermartingale
(Lemma IV.14.3 in \cite{rw-2}).
Since $A_{t\wedge\sigma_A}$ is a decreasing process
(and therefore, a fortiori, a supermartingale),
$e^{-q{(t\wedge\sigma_A})}p(X_{t\wedge\sigma_A})$ is a supermartingale.

We now prove the converse statement.
Assume that $p$ is $q$-superharmonic on $S$.
Hence, for any $A\subset S$ with compact support,
we have \eqref{eq:doob-mayer}, and
$e^{-q{(t\wedge\sigma_A})} p(X_{t\wedge\sigma_A})$ is a supermartingale.
Note that $A_{t\wedge\sigma_A}$ is a predictable
finite variation process and
$M_{t\wedge\sigma_A}$ is a local martingale.
Since $e^{-q{(t\wedge\sigma_A})}p(X_{t\wedge\sigma_A})$ is a càdlàg supermartingale,
there exist a predictable decreasing process
$\widetilde A$
and a local martingale
$\widetilde M$
such that
\[
e^{-q{(t\wedge\sigma_A})}-p(X_{t\wedge\sigma_A})-p(X_0)
=
\widetilde A_t+\widetilde M_t.
\]
(Theorem~ III.16 \cite{protter2012stochastic}.)
Hence, by \eqref{eq:doob-mayer},
\[
A_{t\wedge\sigma_A}-\widetilde A_t
=
\widetilde M_t-M_{t\wedge\sigma_A}
=:N_t.
\]

$N_t$ is a local martingale and a predictable
finite variation process.
In particular, it is a special semimartingale
(Chapter III in \cite{protter2012stochastic}).
$N_t$ admits two different canonical decompositions:
\[
N_t
=
0+N_t
=
N_t+0,
\]
where the first terms are predictable finite variation
processes and the second terms are local martingales.
However, since the decomposition is unique (Theorem~III.34 in \cite{protter2012stochastic}),
$N_t$ must be $0$.
Therefore, the predictable decreasing process $\widetilde{A}$ has the representaiton
\[\widetilde{A}_t = A_{t\wedge \sigma_A} = \int_0^{t\wedge \sigma_A} e^{-qs}(\LL-q)p(X_{s-})\diff s + \frac{1}{2}\sum_{a \in F^p} \Delta p'(a) \int_0^{t\wedge \sigma_A} e^{-qs} \diff L_s^a.\]

Finally, we prove $(\LL-q)p(x) \leq 0$ for $x \in (\R \setminus F^p) \cap A$ and
$\Delta p'(x) \leq 0$ for $x \in F^p \cap A$.
We define the almost surely non-positive Stieltjes measure associated
with $\widetilde A$ by
\begin{align*}
d\widetilde A((t_1,t_2])
:=
\widetilde A_{t_2}-\widetilde A_{t_1}
=&\int_{t_1\wedge\sigma_A}^{t_2\wedge\sigma_A}
e^{-qs}(\mathcal L-q)p(X_{s-})\,ds  \\
&+ \frac{1}{2}\sum_{a \in F^p} \Delta p'(a) \int_{t_1\wedge \sigma_A}^{t_2\wedge \sigma_A} e^{-qs} \diff L_s^a.
\end{align*}

Define the signed measures
\[
m((t_1,t_2])
:=
\int_{t_1\wedge\sigma_A}^{t_2\wedge\sigma_A}
e^{-qs}(\mathcal L-q)p(X_{s-})\,ds,
\]
and, for each $a\in F^p$,
\[
\nu_a((t_1,t_2])
:=
\frac12\Delta p'(a)
\int_{t_1\wedge\sigma_A}^{t_2\wedge\sigma_A}
e^{-qs}\,dL_s^a.
\]
Then
\[
d\widetilde A=m+\sum_{a\in F^p}\nu_a.
\]

Let $x\in F^p\cap A$, and consider a path 
started from $x$.
Since
\(
\{t\leq\sigma_A:X_t=x\}
\)
has Lebesgue measure zero, the absolute continuity of $m$ yields
\(
m\bigl(\{t\leq\sigma_A:X_t=x\}\bigr)=0.
\)
Moreover, for every $a\in F^p$ with $a\neq x$, the defining property of the
local time (Theorem IV. 43. 3 in \cite{rw-2}) implies that $\nu_a$ is concentrated on
\(
\{t\leq\sigma_A:X_t=a\},
\)
and therefore
\(
\nu_a\bigl(\{t\leq\sigma_A:X_t=x\}\bigr)=0.
\)
Consequently,
\[
d\widetilde A\bigl(\{t\leq\sigma_A:X_t=x\}\bigr)
=
\nu_x\bigl(\{t\leq\sigma_A:X_t=x\}\bigr).
\]
Since
\[
\nu_x\bigl(\{t\leq\sigma_A:X_t=x\}\bigr)
=
\frac12\Delta p'(x)
\int_{\{t\leq\sigma_A:X_t=x\}}
e^{-qs}\,dL_s^x,
\]
and
\[
\int_{\{t\leq\sigma_A:X_t=x\}}
e^{-qs}\,dL_s^x>0 \quad \mathbb{P}_x\text{-a.s.},
\]
while $d\widetilde A$ is a non-positive measure, we conclude that
\(
\Delta p'(x)\leq0.
\)

Next, let $x\in A\cap(\R\setminus F^p)$ and consider a path 
started from $x$.
Suppose, to the contrary, that
\(
(\mathcal L-q)p(x)>0.
\)
Since $F^p$ is finite and $(\mathcal L-q)p$ is continuous on
$\R\setminus F^p$, there exists an open neighborhood
$B\subset A$ of $x$ such that
\(
B\cap F^p=\varnothing,\)
\(
(\mathcal L-q)p>0
\)
on $B$.
Since no point of $F^p$ is visited before $\sigma_B$,
\(
\nu_a((0,\sigma_B])=0 \) for $a \in F^p$ $\mathbb{P}_x$-a.s.
Hence
\[
d\widetilde A((0,\sigma_B])
=
m((0,\sigma_B])
=
\int_0^{\sigma_B}
e^{-qs}(\mathcal L-q)p(X_{s-})\,ds
>0,
\]
which contradicts the fact that $d\widetilde A$ is a non-positive measure.
Therefore,
\(
(\mathcal L-q)p(x)\leq0.
\)

Hence, $(\LL-q)p(x) \leq 0$ on $x \in (\R \setminus F^p) \cap A$ and  $\Delta p'(x) \leq 0$ on $x \in F^p \cap A$, which leads to $\mu|_A \geq 0$. Since $A$ is arbitrary, $\mu|_S \geq 0$. 
\end{proof}

Based on this proposition, we define $q$-superharmonicity and $q$-subharmonicity pointwise.
\begin{definition}\label{def:point-super-sub}
The function $p$ is said to be $q$-superharmonic
(resp.\ $q$-subharmonic) at a point $a$ if the following conditions hold:
\begin{enumerate}[label=(\roman*)]
    \item If $a \in \R \setminus F^p$,  then $(\LL-q) p (a) \le 0$
(resp.\ $(\LL-q)  p(a) \ge 0$).
\item  If $a \in  F^p$, then 
$\Delta p(a) \le 0$ (resp. $\Delta p(a) \ge 0$) when $X$ has paths of bounded variation,
and $\Delta p'(a) \le 0$ (resp. $\Delta p'(a) \ge 0$) when $X$ has paths of unbounded variation and $\sigma>0$.
\end{enumerate}
A function is said to be $q$-harmonic if it is both $q$-superharmonic and $q$-subharmonic. It is said to be strictly $q$-subharmonic (resp.\ strictly $q$-superharmonic) if it is $q$-subharmonic (resp.\ $q$-superharmonic) but not $q$-superharmonic (resp.\ $q$-subharmonic).

We refer to $F^p$ as the irregular set of $p$.
\end{definition}

By Proposition \ref{prop:Riesz}, $p$ is $q$-superharmonic on an open set $S$ if and only if it is $q$-superharmonic at every point $a \in S$.
\begin{remark}\label{rem:point-harmonic}
    We observe that $h_a$ is harmonic at $a \in \R$. 
    Let $F$ be the irregular set of $h_a$.
    If $a \notin F$, then
\[
\LL h_a(a)
=
\lim_{a' \downarrow a}\LL h_a(a')
=
0,
\]
where the limit follows from the dominated convergence theorem, since
$h_a$ satisfies Assumption~\ref{assumption-jumpsize}(iv), as shown in
Section~\ref{sec:G-S}.
     If $a \in F$, we have  $\Delta h_a(a)= h_{a}(a+)-g(a-) = 0$ (resp. $\Delta (h_a)'(a)= (h_{a})'(a+)-g'(a-) = 0$) if $X$ has paths of bounded variation (resp. unbounded variation) by using the definition of $h_a$.
\end{remark}

\section{Necessary conditions}\label{sec:necessary}
We present a necessary condition for the optimal stopping problem by using the results in Section.\ref{sec:martin}. 
When $X$ has paths of bounded variation, the continuous fit condition is necessary,
whereas when $X$ has paths of unbounded variation,
the smooth fit condition is required.
This result is consistent with, and can be seen as a generalization of, previous findings such as those by \cite{alili-kyp, Egami-Yamazaki-2011, mordecki}.

Henceforth, we consider only the case where
\begin{align}\label{eq:assuption-q}
    \text{\( q = 0 \) and \( \lim_{t \to \infty} X_t = \infty \) almost surely under \( \p \).}
\end{align}
 This assumption does not entail any loss of generality.  If $q=0$ and $X$ oscillates, 
 the problem is trivial. Since, for any point $x$, the probability of reaching $x$ in finite time is one regardless of the initial distribution, the optimal stopping time $\tau^*$ is given as $\tau^* = \inf \{t \geq 0: X_t = \max g\} $ if $g$ attains its maximum. 
 It is optimal never to stop if $g$ does not attain its maximum.
In cases where  \( q > 0 \) or [$q=0$ and \( \lim_{t \to \infty} X_t = -\infty \)], the problem can be reduced to \( q = 0 \) and \( \lim_{t \to \infty} X_t = \infty \) by the exponential change of measure:
\begin{align}
    \frac{\diff \p^*_x}{\diff \p_x}\Big|_{\F_t}:= e^{\Phi(q)(X_t-x)-qt}. \label{eq:expochange}
\end{align}
Let $\E^*$ denote the expectation operator under $\p^*$.
Then, we have 
\[    \E_x[e^{-qt}g(X_t)] = e^{\Phi(q)x}\E^*_x [g^*(X_t)], \quad  \text{where}\quad  g^*(x)=e^{-\Phi(q)x}g(x)\]
 
The optimal stopping problem for $g$ reduces to that for $g^*$.
A straightforward calculation yields
\begin{align}
\psi^*(\theta) = \psi(\theta+\Phi(q))-q, \conn
\LL^* g^*(x) = e^{-\Phi(q)x} (\LL - q)g(x),  \label{tildeL}
\end{align}
where $\psi^*$ is the Laplace exponent and $\LL^*$ is the generator of $X$ under $\p^*$.
This identity is used in practical applications; see Section \ref{subsection:onesided}.
Moreover, since  $(\psi^*)'(0+) = \psi'(\Phi(q)) >0$ on the account of the strict convexity of $\psi$ and $\Phi(q)>0$, $(X, \p^*)$ always drifts to infinity. Hence this case can be reduced to \( q = 0 \) and \( \lim_{t \to \infty} X_t = \infty \).

\begin{proposition}\label{prop:necessity}
    We assume Assumption \ref{assumption-jumpsize} holds.
     Suppose that $X$ has paths of bounded variation (resp. unbounded variation).
     Let $-\infty<a<b<\infty$.
    If $v(x)=H_a g(x)$ for $x \in (a, \infty)$ and $a \notin F$, then we have
    \begin{align*}
        H_a g (a) = g(a) \quad \Big(\text{resp.} \ (H_a g)'(a+)=g'(a) \Big).
    \end{align*}
    If $v(x)=H_{a,b}g(x)$ for $x \in (a, b)$ and $a, b \notin F$, then we have
    \begin{align*}
    \begin{cases}
         H_{a,b}g (a+) &= g(a) \quad \Big(\text{resp.} \ (H_{a,b}g)'(a+)=g'(a) \Big),\\
         H_{a,b}g (b-) &= g(b) \quad \Big(\text{resp.} \ (H_{a,b}g)'(b-)=g'(b) \Big).
    \end{cases}
    \end{align*}
\end{proposition}
\begin{remark}
If $X$ has paths of unbounded variation, then $H_a g(a) = g(a)$ and  $H_a g(b) = g(b)$ holds automatically, since $a$ and $b$ are regular in the sense that $\p_a[T^{-}_a=0]=1$ and  $\p_b[T_b=0]=1$
    in this case.



\end{remark}

\begin{proof}[Proof of Proposition \ref{prop:necessity}]
    We use the fact that $v$ is superharmonic majorant of $g$ (Theorem 1 in Section 3.3 of  \cite{Shiryaev_2008}). 
     First, we assume that  $v=H_a g$. 
    Consider the case that $X$ has paths of bounded variation (resp. unbounded variation).
    Since $v=H_a g$ is superharmonic on the entire line, the corresponding signed measure $\mu$ in Proposition \ref{prop:Riesz} must be a non-negative measure. Thus, it is necessary that $\mu(\{a\})\ge 0$, which is equivalent to $\Delta H_a  g (a) \le 0$  because we have $\delta >0$ $($resp. $\Delta (H_a  g)' (a) \le 0$ $ \text{due to } \sigma^2 > 0 )$. On the other hand, since $v=H_a g$ is a majorant of $g$ and $g(x)=H_a g(x)$ for $x\leq a$, $\Delta H_a  g (a) \ge 0$ $($resp. $\Delta (H_a  g)' (a) \ge 0$ $)$is necessary. Thus, $\Delta H_a  g (a) = 0$, so that $ H_a g (a) = g(a)$ $($resp. $(H_a g)'(a+)=g'(a)$ $)$ is obtained.
    The second case follows the same argument as the first.
\end{proof}

\section{Sufficient conditions}\label{sec:sufficient}
In this section, we present sufficient conditions for optimality. 
First, we introduce Condition~\eqref{conditionA} and show that, whenever this condition fails, the optimal strategy is to never stop in Proposition \ref{prop:conditionA}.
Next, under Condition~\eqref{conditionA}, Proposition~\ref{prop:verification} establishes the existence of a solution and provides a constructive procedure for obtaining it.

Recall that Assumption 1(ii) requires $g$ to be upper semicontinuous and right-continuous in the bounded variation case, which implies that $\Delta g(x) \geq 0$ for every $x$. Thus, our approach may appear to cover only the case where $\Delta g(x) \geq 0$ for every $x$. However, the proof of Proposition \ref{prop:verification} does not require the upper semicontinuity of $g$. Consequently, in the case where Proposition \ref{prop:verification} applies, the upper semicontinuity assumption is not needed. In fact, Section~\ref{subsection:discontinuous} provides a concrete example in which $g$ is not upper semicontinuous and $\Delta g(a) < 0$ for some $a$.

The key ingredient in the proofs of these propositions is the analysis of superharmonic and subharmonic functions and their analytic properties, including the maximum principle (Proposition \ref{prop:subsuperMP}).

Our solution strategy  consists of successively “eliminating” the strictly subharmonic components starting from $-\infty$.
Let $D$ denote the set of strictly subharmonic points (See Definition \ref{def:point-super-sub});
\[
D = D_R \cup  D_I \quad \text{with} \quad D_R := \{x \notin F: \ \LL g(x)>0 \}
\]
\[ \conn D_I := \{x \in F: \ \Delta g(x)>0 \} \quad (\text{resp.}\  D_I := \{x \in F: \ \Delta g'(x)>0 \}),\]
if $X$ has paths of bounded variation (resp. unbounded variation).
By the continuity of $\LL g$ on $\R \setminus F$, $D_R$ is open. 
We present the following assumptions on the reward function $g$.

\begin{assump}[Assumptions on $g$]\label{assump:g}
\begin{enumerate}[label=(\roman*)]
    \item $D$ has $n < \infty$ connected components.
    \item There exists $\xi$ such that $g$ is strictly subharmonic or strictly superharmonic on $(-\infty, \xi)$.
\end{enumerate}

\end{assump}
Assumption \ref{assump:g}-(ii) is imposed solely to simplify the exposition of Proposition \ref{prop:conditionA}. It is not used elsewhere in the paper. In contrast, Assumption \ref{assump:g}-(i) is more essential, as it guarantees that the iterative solution procedure terminates after finitely many steps.
This assumption is imposed throughout the paper.
Figure \ref{fig:realline2} illustrates a set \(D\) with \(n=4\) connected components.

\input{realline2_modified}




Our approach proceeds by separating the analysis into two cases: whether $g$ attains its maximum or not.
If $g$ attains a maximum,  let $\hat{g} := H_{-\infty, \beta}g$, where
\begin{align}\label{eq:beta}
    \beta:= \sup \{x: g(x) = \max g\ \text{for all } \ y\}.
\end{align}
Observe that
\begin{align}\label{eq:ghat}
        \hat{g}(x)
= g(\beta)\,\mathbf{1}_{(-\infty,\beta]}(x)
  + g(x)\,\mathbf{1}_{(\beta,\infty)}(x).
\end{align}
If $\beta= \infty$, then the optimal strategy is never to stop under the assumption \eqref{eq:assuption-q}.
Henceforth, we assume that $\beta< \infty$.

\begin{figure}[h]\label{fig:flowchart}
  \centering
  \includegraphics[width=0.5\linewidth]{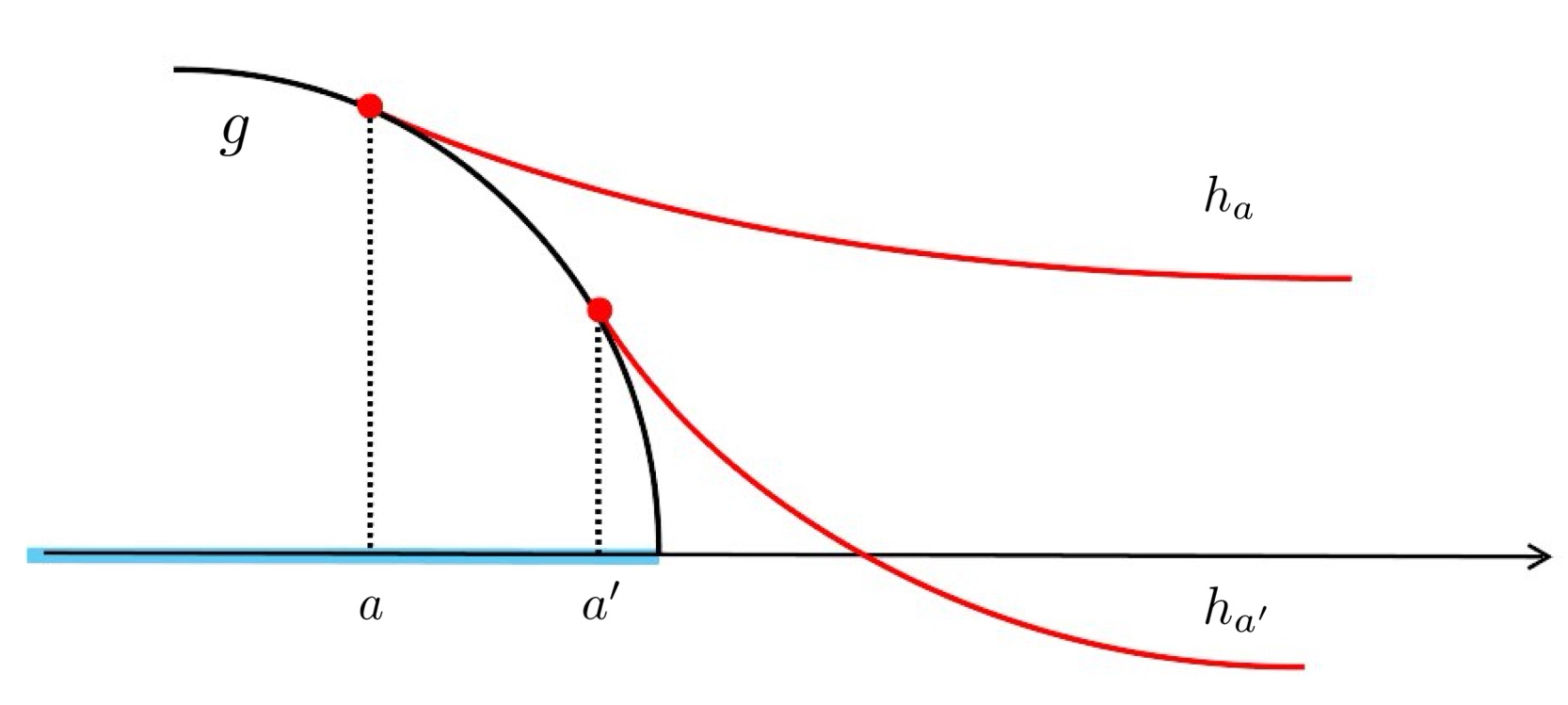}
  \captionsetup{font=scriptsize}
  \caption{Relative position of $h_a$ and $g$ is presented. $h_{a'}$ is not majorant of $g$ while $h_a$ is a majorant of $g$ if we take  $a$ large in the negative direction. }
  \label{fig:functionh_a}
\end{figure}

We consider the following condition on $g$;
\begin{align}
   \inf \{a \in \R: h_a(x)\le g(x)
    \text{ for some }x > a\}= -\infty, \tag{A} \label{conditionA}
\end{align}
where $h_a$ is the smooth Gerber–Shiu function for $q=0$ defined in \eqref{eq:GerberShiu}. 
We refer this condition as Condition \eqref{conditionA} hereafter.
Figure \ref{fig:functionh_a} illustrates the relative positions of $g$ and its smooth Gerber–Shiu function.
Recall the definition of the stopping and continuation region $\Gamma$ and $C$ in \eqref{eq:Gamma-C}.

\begin{proposition}\label{prop:conditionA}
Assumptions \ref{assumption-jumpsize} and \ref{assump:g} are in place.
\begin{enumerate}[label=(\Roman*)]
    \item Assume that $g$ does not attain the maximum.     If Condition \eqref{conditionA} fails, then $C= \R$.
    \item Assume that $g$ attains its maximum. Then, $\hat{g}$ satisfies Condition \eqref{conditionA}.
\end{enumerate}
\end{proposition}

To prove this proposition and those presented in the subsequent sections,
we investigate, in Lemmas~ \ref{lemma:sub-dominance}, \ref{lemma:majorant} and \ref{lemma:monotonicity},
the behavior of the smooth Gerber–Shiu function $h_a$ depending on whether $g$ is superharmonic or subharmonic.

\begin{lemma}\label{lemma:sub-dominance}
     Let $a<b$.
     Suppose that $a \in D$. Then, there exists $\eps>0$ such that $s:=g-h_{a}$ is  strictly positive on $(a,a+\eps)$. 
     Moreover, if $g$ is subharmonic on $(a , b)$, then $s>0$ on $(a,b)$. 
\end{lemma}

\begin{proof}
     Assume that  $X$ has paths of unbounded variation.
     First, we consider the case where $a \in D_I$. 
     Since $a \in D_I$, $g'(a+)>g'(a-)$. By the definition of $h_a$, $g'(a-) = h_a'(a+)$. Therefore, $s'(a+)= g'(a+)-h_a'(a+)>0$. Since $s(a)=0$ by the definition of $h_a$, $s$ is strictly positive on $(a,a+\eps)$ for some $\eps>0$.

    We consider the case where $a \in D_R$. Since $D_R$ is open, we take an interval $(\ell, r) \subset D_R$ with $a \in (\ell, r)$.   
    We have $s =0$ on $(-\infty, a]$, $\LL s >0$ on $(a, r)$ and $s'(a+)=g'(a+)-h_a'(a+)=g'(a+)-g'(a-)=0$
        by $a \in D_R$. 
Take $\delta>0$ with $[a, a+\delta] \subset (\ell, r)$.
We have, for each $x \in (a, a+\delta)$,
    \begin{align*}
        0<c& := \min_{y \in [a, a+\delta]} \LL g(y) \leq  \LL g(x) = \LL s(x) \\
        & =- \gamma s'(x)+ \frac{1}{2}\sigma^2 s''(x)+\int_{(-\infty, 0)}[s(x+y)-s(x)-y\1_{(-1,0)}(y)s'(x)] \Pi(dy),
        \end{align*}
    where the first inequality is obtained from $\LL g >0$ on $(\ell, r)$ and the continuity of $\LL g$ on $[a, a+\delta]$, and the second equality due to the harmonicity of $h_a$.
    Letting $x \downarrow a$, we obtain, from the dominated convergence theorem and the facts that $s'(a+)=0$ and $s =0$ on $(-\infty, a]$,
    \begin{align*}
        0<c &\leq -\gamma s'(a+)+\frac{1}{2}\sigma^2 s''(a+)+\int_{(-\infty, 0)}[s(a+y)-s(a)-y\1_{(-1,0)}(y)s'(a+)] \Pi(dy)\\
        &=\frac{1}{2}\sigma^2 s''(a+),
    \end{align*}
    and thus $s''(a+)>0$. This and $s'(a+)=0$ imply that $s>0$ on $(a, a + \eps)$ for some $\eps>0$.
    Finally, assume that $g$ is subharmonic on $(a,b)$. 
    $s$ is also subharmonic on $(a,b)$.
    By Proposition \ref{prop:subsuperMP}, for any $d \in (a,b)$, $s|_{(-\infty, d)}$ does not attain its maximum; therefore, 
    $s$ is  strictly increasing on $(a,b)$. In particular,  $s>0$ on $(a,b)$. 

Next, we consider the bounded variation case. The proof proceeds similarly to the unbounded variation case, with the main differences being the order of the derivatives evaluated and the leading term in the limit of the operator $\mathcal{L}$.
First, for the case where $a \in D_I$, we have $g(a+)>g(a-)$. By the definition of $h_a$, this yields $s(a+)= g(a+)-h_a(a+)>0$. The right-continuity of $s$ then guarantees that $s>0$ on $(a,a+\eps)$ for some $\eps>0$.
Second, we consider the case where $a \in D_R$. In this case, we have $s(a+)=0$. Following the same procedure as above, we evaluate $\mathcal{L} s(x)$ for $x \in (a, a+\delta)$ and let $x \downarrow a$. By the dominated convergence theorem and the fact that $s=0$ on $(-\infty, a]$, we obtain

\begin{align*}    
0<c \leq \delta s'(a+)+\int_{(-\infty, 0)}[s(a+y)-s(a)] \Pi(dy)= \delta s'(a+). 
\end{align*}
Thus, we have $s'(a+)>0$, which, combined with $s(a+)=0$, implies $s>0$ on $(a, a + \eps)$ for some $\eps>0$.
\end{proof}

\begin{lemma}\label{lemma:majorant}
Let $a<b$. Suppose that $g$ is strictly superharmonic on $a$.
Then, there exists $\eps>0$ such that $f := h_{a}-g$ is  strictly positive on $(a,a+\eps)$.
Moreover, if $g$ is superharmonic on $(a , b)$, then $f>0$ on $(a,b)$. 
\end{lemma}
\begin{proof}
The proof is exactly the same as that of Lemma \ref{lemma:sub-dominance}.    
\end{proof}


\begin{lemma}\label{lemma:monotonicity}    Let $a<b$.  Suppose that $g$ is strictly superharmonic on $a$ and superharmonic on $(a,b)$. 
Then,  we have $h_{a}>h_{b}$ on $(a, \infty]$.
In particular, for each $x>a$, $h_{a'}(x)$ converges increasingly to $h_{a}(x)$ as $a'\downarrow a$.
\end{lemma}
\begin{proof}     
It follows from Lemma \ref{lemma:majorant} that $f := h_{a} - g$ satisfies $f = 0$ on $(-\infty, a]$ and $f > 0$ on $(a, b)$.     
Define $h := h_{a} - h_{b}$. Then $h = 0$ on $(-\infty, a]$ and $h = f > 0$ on $(a, b)$.
We show that $h$ is subharmonic on $(a,\infty)$.
By assumption, $h$ is subharmonic on $(a,b)$.
Moreover, since both $h_a$ and $h_b$ are harmonic on $(b,\infty)$, so is $h$.
Thus, it remains to prove that $h$ is subharmonic at $b$.
$h_a$ is harmonic on $b$ by the definition of the smooth Gerber-Shiu function.
Moreover $h_b$ is also harmonic on $b$  from Remark \ref{rem:point-harmonic}.
Therefore, $h$ is harmonic at $b.$

Since $h$ is subharmonic on $(a, \infty)$, $h$ restricted on $(-\infty, d)$ does not attain its maximum for any $d>a$ by the maximum principle (Proposition \ref{prop:subsuperMP}).
Hence, $h$ must be strictly positive on $(b, \infty]$.
\end{proof}

\begin{proof}[Proof of Proposition \ref{prop:conditionA}]
(I)\  Suppose that Condition \eqref{conditionA} fails. 
Under Assumption \ref{assump:g}, for sufficiently negatively large $\xi$, $(-\infty, \xi)$ is (i) strictly subharmonic or (ii) strictly superharmonic.

First, we show that there exists $r$ such that 
\begin{align}
  (-\infty, r) \subset C.  \label{eq:minus-inf-r}
\end{align}
Consider the case (i).
Take $x \in (-\infty, \xi)$.
There exists an open interval $(x-\eps, x+\eps)$ with $(x-\eps, x+\eps) \subset (-\infty, \xi)$.
Since $g$ is strictly superharmonic on $(x-\eps, x+\eps)$, we have
\[v(x) \geq H_{x-\eps, x+\eps} g(x) > g(x), \]
which leads to $x \in C$.
Hence, $(-\infty, \xi) \subset C$.
Therefore, we have $r$ defined in \eqref{eq:minus-inf-r}.

Consider the case (ii).
By the assumption that Condition \eqref{conditionA} fails and Lemma \ref{lemma:monotonicity},
there exists $r<\xi$ such that for any $r' < r$  
\begin{align}\label{eq:xi-prime}
     h_{r'}(y) \leq g(y) \quad \text{for some}\  y> r'.
\end{align}
We have, for $x \in (r', r)$,
\begin{align*}
    v(x) \geq H_{r', y} g(x)
    = h_{r'}(x)+W(x-r')\frac{g(y)-h_{r'}(y)}{W(y-r')}
    \geq h_{r'}(x)
    > g(x),
\end{align*}
where the second equality follows from \eqref{eq:expectedreward2}, the third inequality follows from \eqref{eq:xi-prime} and the final inequality follows from Lemma \ref{lemma:majorant}. 
Hence, $(r', r) \subset C$.
Since $r'<r$ is arbitrary, we obtain $(-\infty, r) \subset C$.

It follows that there exists $r_0 >r$ such that $(-\infty, r_0)$ is a connected component of $C$.
Since we assume that $g$ is upper semi-continuous, $v = H_\Gamma g$ by Theorem 3.3.3 in \cite{Shiryaev_2008}. 
 Since the process starting from $x \in (-\infty, r_0)$ exits $C$ from $r_0$ in finite time  from the assumption $X_t \goes$ as $t \goes$ and the support of $\p_x [X_{T_{(-\infty,r_0)^\mathrm{c}}} \in \diff y]$ is $\{r_0\}$ (recall $X$ is spectrally negative), we obtain, for $x \in (-\infty, r_0)$,
\[ g(x) \leq v(x) = H_\Gamma g(x) = H_{(-\infty,r_0)^\mathrm{c}} g(x) = g(r_0).\]
Since we assume that $g$ does not attain its maximum, $\limsup_{x \goes} g(x)=\infty$.

Take $z \in \R$ arbitrarily.
There exists $y$ such that $g(z)<g(y)$.
Hence, since we have $g(z) <g(y) = H_{-\infty, y} g(z) \leq v(z)$, it holds $z \in C$.
Since $z$ is taken arbitrarily, we conclude $C= \R$.

(II)\   $\hat{g}$ is (super)hamonic on $(-\infty, \beta)$ since it is a constant on $(-\infty, \beta)$.
 Since $\hat{g}(x) = \hat{g}(\beta)$ for $x \in (-\infty, \beta]$ and $g'(\beta-)=0$, 
 it follows from the definition of the Gerber-Shiu function in \eqref{eq:GerberShiu} that,
 for $x>\beta$, $\hat{h}_\beta(x) =\hat{g}(\beta)>\hat{g}(x)$, where $\hat{h}_\beta$ is the smooth Gerber–Shiu function for $\hat{g}$ at $\beta$. Therefore, Condition \eqref{conditionA} holds.
\end{proof}

\subsection{Sufficient Condition}

\begin{proposition}\label{prop:verification}
Assume that $g$ does not attain its maximum and that Condition~\eqref{conditionA} is satisfied.
Set $g^0=g$. We recursively define, for $i\ge1$,
\begin{align}
a_i &:= \inf \{a \in \R: h^{i-1}_a(x)\le g^{i-1}(x)
    \text{ for some }x > a\}, \label{eq:a_i}\\
\alpha_i &:= \max\{\alpha\in\R:
h^{i-1}_{a_i}(x)-\alpha W(x-a_i)\ge g^{i-1}(x)
\text{ for all }x >  a_i\}, \label{eq:alpha_i}\\
a'_i &= \inf \{a>a_i:h^{i-1}_{a_i}(a)-\alpha_iW(a-a_i)\neq g^{i-1}(a) \}\\
b_i &:= \sup\{b>a'_i:
h^{i-1}_{a_i}(b)-\alpha_iW(b-a_i)=g^{i-1}(b)\},\label{eq:b_i}\\
g^i &:=  H_{a_i,b_i}g^{i-1} \label{eq:g^i},
\end{align}
as long as both $a_i$ and $b_i$ are finite, where $h^i_a$ is the smooth Gerber-Shiu function of $g^i$ at $a$.
Let
\(
m=\min\{i\ge1:\ a_i=\infty \text{ or } b_i=\infty\}
\) and $E = \bigcup_{1 \leq i \leq m} [a_i, b_i]$,
where we set $[a_m, b_m] = \varnothing$ if $a_m = \infty$.
Then, 
\[ v = g^m =  H_{\R \setminus E} g \conn \tau^*=T_{\R \setminus E}, \]
and we have $m\le n<\infty$, where $n$ is the number of subharmonic components.
\end{proposition}

We define, for each $i \geq 1$, 
\[D_R^{i-1} := \{x \notin F^{i-1}: \ \LL g^{i-1}(x)>0 \} \conn D_I^{i-1} := \{x \in F^{i-1}: \ \Delta g^{i-1}(x)>0 \}, \]
where $F^{i-1}$ is the finite set of all point where the first or second derivative of $g^{i-1}$ is not defined. 
We denote $D^{i-1} = D_R^{i-1} \cup D_I^{i-1}$.
In particular, $D^0 = D$, $D_R^0=D_R$, $D_I^0=D_I$ and $F^0 = F$.

We prove Proposition \ref{prop:verification} through a sequence of lemmas.
To avoid repetition, we present the proofs in these lemmas only for the case of unbounded variation
The parallel arguments for the bounded variation case follow the same pattern as demonstrated in Lemma \ref{lemma:sub-dominance}.

\begin{lemma}\label{lemma:a_i-b_i}
    For each $i = \{ 1, \cdots , m-1\}$, $a_i, b_i \in \R \setminus D^{i-1}$ and $(b_i, a_{i+1}) \cap D^{i} = \varnothing$.
\end{lemma}
\begin{proof}
    Since the argument for $i\geq 2$ is identical, we only prove the case $i=1$. We first show that $a_1 \in \R \setminus D$. Assume to the contrary that $a_1 \in D_R$ or $a_1 \in D_I$. First, we consider the case where $a_1 \in D_R$. Since $D_R$ is open, there exists $\eps>0$ such that $a_1 -\eps \in D_R$. It follows from Lemma \ref{lemma:sub-dominance} that $h_{a_1-\eps}(x)<g(x)$ for some $x > a_1-\eps$, which contradicts the definition of $a_1$. We consider the case where $a_1 \in D_I$. By Lemma  \ref{lemma:sub-dominance}, $h_{a_1}(x)<g(x)$ for some $x > a_1$. Since the mapping $a \mapsto h_a$ is left-continuous, $h_{a_1-\eps}(x)<g(x)$ for sufficiently small $\eps>0$, which also contradicts the definition of $a_1$.

    We prove that $b_1 \in \R \setminus D$.
    By the definition of $\alpha_1$, we have, for each $x>a_1$,
    \begin{align}
        f^1(x) := h_{a_1}(x)-\alpha_1 W(x-a_1) \geq g(x). \label{eq:f1>=g}
    \end{align}
    Moreover,  $f^1(b_1) =  g(b_1)$ by the definition of $b_1$.
    Since both $h_{a_1}$ and $W$ are differentiable at $b_1$ (see  Sections \ref{sec:scale_functions} and \ref{sec:G-S}), so is $f^1$.
    Hence, 
    \begin{align}
      g'(b_1-) \geq (f^1)'(b_1) \geq g'(b_1+). \label{eq:smoothfit-b1}
    \end{align}
    Assume that $b_1 \in F$.
    Since $g'(b_1-) \geq g'(b_1+)$ from \eqref{eq:smoothfit-b1}, it follows from Proposition \ref{prop:Riesz} that $b_1 \in \R \setminus D$. 
    Then, assume that $b_1 \notin F$
    By \eqref{eq:smoothfit-b1}, $g'(b_1)= (f^1)'(b_1)$.
    Combining this with \eqref{eq:f1>=g},
    we obtain $g''(b_1) \leq  (f^1)''(b_1)$. Hence, using \eqref{eq:operator-L},
   $\LL g(b_1) \leq \LL f^1 (b_1) = 0.$
   Therefore, $b_1 \in \R \setminus D$.
    
    Finally, we prove $(b_1, a_2) \cap D^1 = \varnothing$. 
    It follows from Lemma \ref{lemma:sub-dominance} that, for each $x \in D^1 \cap (b_1, \infty)$, $h_x < g^1$ for some points. Therefore, $x \geq a_2$; and hence, $(b_1, a_2) \cap D^1 = \varnothing$.
\end{proof}

\begin{lemma}\label{lemma:expectedreward}
     For each $i = \{ 1, \cdots , m-1\}$ and  $x \in (a_i, b_i)$, we have $g^m(x) = g^i(x) =h^{i-1}_{a_i}(x)-\alpha_i W(x-a_i).$
\end{lemma}
\begin{proof}
    The first equality directly follows from the definitions of $g^i$ and $g^m$.
    For the second equality, we have 
    \begin{align*}
       g^i(x) &=  H_{a_i, b_i} g^{i-1}(x) = h^{i-1}_{a_i}(x) - W(x-a_i) \frac{g^{i-1}(b_i)-h^{i-1}_{a_i}(b_i)}{W(b_i-a_i)}\\
       &=  h^{i-1}_{a_i}(x) - \alpha_i W(x-a_i),
    \end{align*}
    where the second equality follows from the definition of $b_i$ in \eqref{eq:b_i} and \eqref{eq:expectedreward2}. To obtain the third equality, we use the following identity followed by \eqref{eq:b_i} and \eqref{eq:alpha_i};
    \begin{align}
        h^{i-1}_{a_i}(b_i)-\alpha_i W(b_i-a_i)=g^{i-1}(b_i). \label{eq:alpha-identity}
    \end{align}  
\end{proof}

\begin{lemma}
    $g^m = H_{\R \setminus E} g$. \label{lemma:DP}
\end{lemma}
\begin{proof}
    Let $E_i = \bigcup_{j=1}^i (a_j, b_j)$.
    It suffices to show that
    \begin{align}
        g^i(x) = H_{\R \setminus E_i} g(x) \label{eq:math-ind}
    \end{align}
    holds for all $i \in \{1, \dots, m\}$.
    We proceed by mathematical induction on $i$.
    
    For $i=1$, \eqref{eq:math-ind} holds by the definition of $g^1$ in \eqref{eq:g^i}.
    Assume that \eqref{eq:math-ind} holds for a fixed $i \le m-1$. 
    Then, for every $x \in \R$,
    \begin{align*}
        g^{i+1}(x) &= H_{a_{i+1}, b_{i+1}} g^i(x) = H_{a_{i+1}, b_{i+1}} H_{\R \setminus E_i} g(x)\\
        &= \int \int g(z) \p^y[X_{T_{\R \setminus E_i}} \in \diff z] \p^x[X_{T_{\R \setminus (a_{i+1}, b_{i+1})}} \in \diff y].
    \end{align*}
    Let $A := (-\infty, a_{i+1}] \cup [b_{i+1}, \infty)$.
    Since $X$ is spectrally negative, the exit distribution $\p^x[X_{T_{\R \setminus (a_{i+1}, b_{i+1})}} \in \diff y]$ is supported on $A$. 
    Moreover, since $(\R \setminus E_i) \cap A = (\R \setminus E_{i+1}) \cap A$, the inner transition probability satisfies
    \[
        \p^y[X_{T_{\R \setminus E_i}} \in \diff z] = \p^y[X_{T_{\R \setminus E_{i+1}}} \in \diff z] \quad \text{for all } y \in A.
    \]
    Substituting this into the integral and applying the strong Markov property along with the inclusion $\R \setminus E_{i+1} \subset \R \setminus (a_{i+1}, b_{i+1})$, we obtain
    \begin{align*}
        g^{i+1}(x) &= \int \int g(z) \p^y[X_{T_{\R \setminus E_{i+1}}} \in \diff z] \p^x[X_{T_{\R \setminus (a_{i+1}, b_{i+1})}} \in \diff y]\\
        &= \int g(z) \p^x[X_{T_{\R \setminus E_{i+1}}} \in \diff z] = H_{\R \setminus E_{i+1}} g(x).
    \end{align*}
    This completes the induction step and concludes the proof.
\end{proof}

\begin{lemma}\label{lemma:majorant-veri}
    $g^{i+1}$ dominates $g^i$ for each $i$. In particular, $g^m$ is a majorant of $g$.
\end{lemma}
\begin{proof}
    We only show that $g^1$ dominates $g$. The remaining can be shown by the identical argument.
    Since  $g=g^1$ on $\R \setminus (a_1, b_1)$ by  the definition of $g^1$, we prove $g^1 \geq g$ on $(a_1, b_1)$.
    We have, for $x \in (a_1, b_1)$,
    \begin{align*}
        g^1(x) = h_{a_1}(x)-\alpha_1 W(x-a_1) \geq g(x),
    \end{align*}
    where the first equality follows from  Lemma \ref{lemma:expectedreward}, and the second inequality follows from \eqref{eq:f1>=g}.
\end{proof}

\begin{lemma}\label{lemma:Dij}
    For each $i$ and $j$, let $J := (-\infty, a_{\min\{i,j\}+1})$. Then, we have $D^i \cap J=D^j \cap J.$
\end{lemma}
\begin{proof}
    This statement directly follows from the fact that $g^i=g^j$ on $J$. 
\end{proof}

\begin{lemma}\label{lemma:exessive}
    $\{ g^m(X_t): \ t \geq 0 \}$ is a right-continuous supermartingale.
\end{lemma}
\begin{proof}
The right continuity follows from Assumptions \ref{assumption-jumpsize} (ii) and (iii); see the comment following Assumption \ref{assumption-jumpsize}.
It remains to show that $g(X)$ is a supermartingale, which follows if $g$ is superharmonic on the entire real line.
From Proposition \ref{prop:Riesz} with $p=g^m$, it is sufficient to show that $\LL g^m(x) \leq 0$ for $x \in \R \setminus F^m$ and $\Delta (g^m)' (x) \leq 0$ for $x \in F^m$, which is equivalent to $ x \in \R \setminus D^m$.
Let $x \in (-\infty, a_1)$. By the definition of $a_1$ in \eqref{eq:a_i}, $h_x >g$ on $(x, \infty)$. It follows from Lemma \ref{lemma:sub-dominance} that $x \in (\R \setminus D) \cap (-\infty, a_1)$.
By Lemma \ref{lemma:Dij}, we obtain $x \in (\R \setminus D^m) \cap (-\infty, a_1) \subset \R \setminus D^m$.

We next consider the case $x=a_1$.
Suppose first that $a_1 \in F^m$.
It follows from the left-continuity of $a \mapsto h_a(x)$  and the definition of $a_1$ in \eqref{eq:a_i} that $h_{a_1}$ dominates $g$, which implies that $\alpha_1 \geq 0$. 
Then
\begin{align*}
(g^m)'(a_1-)&=g'(a_1-)=h'_{a_1}(a_1+)
\ge h'_{a_1}(a_1+)-\alpha_1W'(0+)\\
&=(g^1)'(a_1+)
=(g^m)'(a_1+),    
\end{align*}

and hence $\Delta (g^m)'(a_1)\le0$.
Now suppose that $a_1\in\mathbb{R}\setminus F^m$. 
Then,
\begin{align}\label{eq:LL-a_1}
    \mathcal{L}g^m(a_1)
=\lim_{z\uparrow a_1}\mathcal{L}g^m(z)
\le0,
\end{align}
where the equality follows from the continuity of $\mathcal{L}g^m$, and the inequality follows from the argument for $(-\infty,a_1)$.

Let $x\in(a_1,b_1)$. By Lemma~\ref{lemma:expectedreward},
\[
g^m(x)=h_{a_1}(x)-\alpha_1W(x-a_1),
\]
which is harmonic. Therefore,
\(\mathcal{L}g^m(x)=0.\)

Next, we consider the case $x = b_1$. 
Suppose first that $b_1 \in F^m$. Then, by using Lemma \ref{lemma:expectedreward} and \eqref{eq:smoothfit-b1}
\[(g^m)'(b_1-)
= (h_{a_1})'(b_1-)-\alpha_1 W'(b_1-a_1)
\geq g'(b_1+) = (g^m)'(b_1+), \]
where the final equality follows from $g^m = g$ on $(b_1, a_2)$. Therefore, we obtain $\Delta (g^m)'(b_1) \leq 0$.
Now suppose that $b_1\in\mathbb{R}\setminus F^m$. Then, as in \eqref{eq:LL-a_1}, we have 
\[ \LL g^m(b_1) = \lim_{z \uparrow b_1} \LL g^m(z) = 0,\]

Let $x \in (b_1, a_2)$. By Lemmas \ref{lemma:a_i-b_i}, $x \in \R \setminus D^1$. Morerover, since $D^1 \cap (b_1, a_2)= D^m \cap (b_1, a_2)$ from Lemma \ref{lemma:Dij}, we obtain $x \in  \R \setminus D^m$.
By repeating the same argument for $x \geq a_2$, we can prove $\LL g^m(x) \leq 0$ for $x \in \R \setminus F^m$ and $\Delta g^m (x) \leq 0$ for $x \in F^m$.
\end{proof}

\begin{lemma}\label{lemma:value}
    $g^m$ is the value function and the optimal stopping time is $\tau^* = T_{\R \setminus E}$ .
\end{lemma}
\begin{proof}
    By Lemmas~\ref{lemma:DP}, \ref{lemma:majorant-veri} and \ref{lemma:exessive}, we can apply Lemma~11.1 of \cite{Kyprianou_2014}, which implies that $g^m$ is the value function and that $T_{\R \setminus E}$ is the optimal stopping time.
\end{proof}

\begin{proof}[Proof of Proposition \ref{prop:verification}]
It remains to prove that $m \le n$.
Suppose, to the contrary, that $m>n$.
$g^m$ is superharmonic on $\R \setminus E$ from Lemma \ref{lemma:exessive}.
Since $g=g^m$ on $\R \setminus E$ by the construction of $g^m$ and $g \leq g^m$ by Lemma \ref{lemma:majorant-veri},
$g$ is also superharmonic on $\R \setminus E$. 
Hence, we have $D\subset E$, 
which implies that  each interval $[a_i,b_i]$ contains at least one connected component of $D$.
Hence, there exists $i\in\{1,\ldots,m\}$ such that $[a_i,b_i]$ is disjoint from $D$.
Therefore, for $x\in(a_i,b_i)$,
\[
v(x)=H_Eg(x)= H_{a_i,b_i}g(x)\le g(x).
\]
This contradicts the construction of $a_i$ and $b_i$.
\end{proof}


\begin{remark}
(i): \  At first glance, the precise role of $\alpha_i$ in \eqref{eq:alpha_i} may not be immediately obvious.
          However, we choose $\alpha_i$ satisfying \eqref{eq:alpha-identity} so that $g^m$ can be properly linked to the expected reward $H_{a_i, b_i} g$ on $[a_i, b_i]$ (see Lemma \ref{lemma:expectedreward}).

(ii): \ As seen in the proof of Lemma \ref{lemma:exessive}, $\alpha_i \geq 0$.
    Assume that $X$ has paths of bounded variation (resp. unbounded variation).Since $h_{a_i}^{i-1}(a_i+) = g(a_i-)$ (resp. $(h_{a_i}^{i-1})'(a_i+) = g'(a_i-)$) by the definition of the smooth Gerber-Shiu function, it follows from \eqref{eq:expectedreward2} and \eqref{eq:alpha-identity} that
    \begin{align*}
        H_{a_i, b_i} g(a_i+) - g(a_i-) = -\alpha_i W(0+) \\
    \big(\text{resp.} \ (H_{a_i, b_i} g)'(a_i+) - g'(a_i-) = -\alpha_i W'(0+).\big)
    \end{align*}
    Hence, $H_{a_i, b_i} g$ and $g$ continuously (resp. smoothly) fit at $a_i$ if $\alpha_i = 0$.
    Conversely, if $\alpha_i > 0$, then they do not continuously (resp. smoothly) fit at $a_i$.

(iii):\  At first glance, it might seem natural to define $b_i$ using $a_i$ directly:
    $$b_i:=\sup\{b>a_i : h^{i-1}_{a_i}(b)-\alpha_i W(b-a_i)=g^{i-1}(b)\}.$$
    However, this definition can fail if $g^{i-1}$ is harmonic on some interval $(a_i, a_i+\eps)$ for $\eps>0$. In this case, $b_i$ would incorrectly collapse to $a_i$ because $h^{i-1}_{a_i}$ coincides with $g^{i-1}$ on $(a_i, a_i+\eps)$.
    To avoid this issue, we replace $a_i$ with $a'_i$ in \eqref{eq:b_i}, where $a'_i$ is defined as the right endpoint of the maximal interval to the right of $a_i$ on which $g^{i-1}$ remains harmonic.   
\end{remark}

In Proposition \ref{prop:verification}, we assumed that $g$ does not attain its maximum. When $g$ attains a maximum, the solution can be constructed by applying Proposition \ref{prop:verification} to $\hat{g}$. We summarize this result in the following corollary.

\begin{corollary}
    \label{coro:verification-max}
Assume that $g$ attains its maximum, $\beta$ defined in \eqref{eq:beta} is finite and that Condition~\eqref{conditionA} is satisfied.
Set $\hat{g}^0=\hat{g}$. We recursively define, for $i\ge1$,
\begin{align*}
\hat{a}_i &:= \inf \{a \in \R: \hat{h}^{i-1}_a(x)\le \hat{g}^{i-1}(x)
    \text{ for some }x > a\}, \\
\hat{\alpha}_i &:= \max\{\alpha\in\R:
\hat{h}^{i-1}_{a_i}(x)-\alpha W(x-a_i)\ge g^{i-1}(x)
\text{ for all }x >  a_i\}, \\
\hat{a}'_i &= \inf \{a>\hat{a}^{i-1}_i:\hat{h}_{\hat{a}_i}^{i-1}(a)-\hat{\alpha_i}W(a-\hat{a}_i)\neq \hat{g}^{i-1}(a) \}\\
\hat{b}_i &:= \sup\{b>\hat{a}'_i:
\hat{h}_{\hat{a}_i}^{i-1}(b)-\hat{\alpha}_iW(b-a_i)=\hat{g}^{i-1}(b)\},\\
\hat{g}^i &:=  H_{a_i,b_i}\hat{g}^{i-1} ,
\end{align*}
as long as both $\hat{a}_i$ and $\hat{b}_i$ are finite, where $\hat{h}^i_a$ is the smooth Gerber-Shiu function of $\hat{g}^i$ at $a$.
Let
\(
m=\inf\{i\ge1:\ \hat{a}_i=\infty \text{ or } \hat{b}_i=\infty\}
\) and $\hat{E} = (-\infty, \beta] \cup \bigcup_{1 \leq i \leq m} [\hat{a}_i, \hat{b}_i]$,
where we set $[\hat{a}_m, \hat{b}_m] = \varnothing$ if $\hat{a}_m = \infty$.
Then, 
\[ v = \hat{g}^m =  H_{\R \setminus \hat{E}} g \conn \tau^*=T_{\R \setminus \hat{E}}, \]
and we have $m\le n<\infty$, where $n$ is the number of subharmonic components.
\end{corollary}
\begin{proof}
   The desired results follow by applying Proposition \ref{prop:verification} to $\hat{g}$.
\end{proof}

If there exists a boundary point $\ell$ such that the function is strictly superharmonic on the left of $\ell$ and strictly subharmonic on the right of $\ell$, then the iteration in Proposition~\ref{prop:verification} terminates after one step, i.e., $m=1$.
Moreover, the optimal boundary $a_1$ can be characterized as a solution to an equation, provided that $a_1\notin F$.
Although the range of applicability of this equation is limited, it provides a useful analytical tool for solving the problem.
We use this approach in Examples~\ref{subsection:onesided} and \ref{subsection:discontinuous}.
 We summarize this result in the following corollary.

\begin{corollary}\label{coro:one-sided}
Suppose that $X$ has paths of bounded variation (resp.\ unbounded variation).
Assume that, for some fixed $i$, 
there exists $\ell$ such that $g$ is strictly superharmonic on $(-\infty, \ell)$ and strictly subharmonic on $(\ell, \infty)$.
Assume that the limit
\(M:=\lim_{x\to\infty}g(x)\) exists and is finite.
Let $\kappa$ denote the $\kappa$-function associated with $g^i$ (see \eqref{eq:kappa} and the discussion following it).
Then \(b_{i+1}=\infty\)
(and hence $m=i+1$), and
\[
a_{i+1}=\inf\{a:\kappa(a)\le M \psi'(0) \}.
\]
In particular, if $g$ is continuous (resp. differentiable) at $a_{i+1}$, then $a_{i+1}$ is the unique root of $\kappa-M\psi'(0)$.
\end{corollary}

\begin{proof}
   To simplify the notation, we consider only the case $i=0$ and let $\kappa$ denote the $\kappa$-function associated with $g$. The general case follows from the same argument.
   Let 
   \[ \eta := \inf\{a:\kappa(a)\le M \psi'(0)\}. \]
     Our first statement is that $a_1 =  \eta$.
   First, we prove $a_1 \geq  \eta$.
   Take any $a$ such that $h_a(x) \leq g(x)$ for some $x>a$.
   Let $x^* = \inf \{x: \ h_a(x) \leq g(x)\}$.
   It follows from Lemma \ref{lemma:majorant} that $x^* \in D$ and $x^* \geq \ell$.
   Let $f:=h_a -g$ and assume to the contrary that $\lim_{x \goes} f(x) >0$.
   Then, $f$ attains its minimum on $\{x: \ x \geq x^* \} \subset D$,
   which contradicts that $f$ is superharmonic on $D$ by the maximum principle (Proposition \ref{prop:subsuperMP}).
   Therefore, $\lim_{x \goes} h_a(x) \leq \lim_{x \goes} g(x) = M $, which is equivalent to $\kappa(a)\le M \psi'(0)$ by \eqref{eq:kappa} and \eqref{eq:limit-W} with $q=0$. 
   Hence,
\[
\{a: h_a(x)\le g(x)\ \text{for some }x>a\}
\subset
\{a:\kappa(a)\le M\psi'(0)\}.
\]
Taking the infimum of both sets, we obtain \(a_1\ge \eta\).

   We then prove  that $a_1 \le  \eta$. 
   Assume to the contrary that $a_1 >  \eta$.
   There exists $\delta>0$ such that $\kappa(a_1- \delta) \leq M \psi'(0)$, which is equivalent to 
   \[\lim_{x \goes} h_{a_1-\delta}(x) \leq M.\]
   It follows from Lemma \ref{lemma:a_i-b_i} that $a_i \leq \ell$. 
   Hence, $g$ is strictly superharmonic on $(a_1-\delta, a_1-\delta/2)$.
   Applying Lemma \ref{lemma:monotonicity} with $(a,b)= (a_1-\delta, a_1-\delta/2)$,
   \[ \lim_{x \goes} h_{a_1-\delta/2}(x) < \lim_{x \goes} h_{a_1-\delta}(x)  \leq  M = \lim_{x \goes} g(x).\]
   Hence, $h_{a_1-\delta/2}(x) < g(x)$ for some $x$, which contradicts the definition of $a_1$.
   We conclude $a_1 =  \eta.$

   Next, we prove that $b_1=\infty$. Assume, to the contrary, that $b_1<\infty$. By Lemma \ref{lemma:a_i-b_i}, we have $b_1 \in \mathbb{R}\setminus D$ and $b_1 \le \ell$. On the other hand, Lemma \ref{lemma:monotonicity} yields $h_{a_1}>g$ on $(a_1,\ell)$, which implies that $\ell<b_1$, a contradiction. Therefore, $b_1=\infty$.
   The final part follows from the facts that $\kappa$ is continuous on $\R \setminus F$.
\end{proof}

Finally, we discuss the relationship between our approach and the log-concavity of the reward function. 
Representative works focusing on this property are \cite{hsiau2014logconcave, lin2019one}.
    
\begin{corollary}\label{coro:logconcave}
 Assumptions \ref{assumption-jumpsize} and
 \ref{assump:g} are in force.
   Suppose that  $g$ is log-concave and non-increasing.
   Assume further that $g$ is either strictly log-concave or strictly decreasing.
   Let $x_0 := \sup \{x \in \R: g(x)>0 \}$.
   If $x_0<\infty$, the value function and the optimal stopping time are $v=H_{a_1} g$ and $\tau^* = T_{(-\infty, a_1)}$, where $a_1 := \inf \{a; \kappa(a) \leq 0 \}.$
   If $x_0 = \infty$, then $v=g$ and the optimal stopping rule is to stop immediately.
\end{corollary}
\begin{proof}[Proof of Proposition \ref{coro:logconcave}]
      We prove the result under the assumption  that $g$ is strictly log-concave.
      In the case where $g$ is strictly decreasing, we can prove in the same way.
      We have, for $a < x_0$, $\kappa(a) = g(a) P(a)$, where
\begin{align*}
     P(a) := \Big[\frac{\sigma^2}{2}\frac{g'(a)}{g(a)}+\psi'(0+)-\int_0^\infty \diff x \int_{(x, \infty)}\Pi (-\diff z)\Big(\frac{g(x+a-z)}{g(a)}-1\Big)\Big].
\end{align*}
It follows from the strict log-concavity of $g$ and a direct calculation that $g'(\cdot)/g(\cdot)$ is strictly decreasing and for every $z>x$, $g(x+\cdot-z)/g(\cdot)$ is strictly increasing.
Hence, $P(a)$ is strictly decreasing on $(-\infty, x_0)$.
Since $g$ is non-increasing by the assumption, $\kappa$ is strictly decreasing on $(-\infty, x_0)$.
By Lemma \ref{lemma:kappa-dec} stated immediately after this proof, $g$ is strictly superharmonic  on $(-\infty, x_0)$ .
Therefore, Corollary \ref{coro:one-sided} yields the desired result if $x_0 < \infty$.  
If $x_0= \infty$, then $g$ is superharmonic on $\R$ and $v=g$.
\end{proof}

\begin{lemma}\label{lemma:kappa-dec}
    $g$ is superharmonic  on $(a_1, a_2)$ if and only if $\kappa$ is decreasing  on $(a_1, a_2)$.
\end{lemma}
\begin{proof}
    We only prove the unbounded variation case.
    Assume that $g$ is superharmonic on $(a_1, a_2)$.
   For $a \in (a_1, a_2) \setminus F$, we have $\LL g(x) \leq 0$. It follows from \eqref{eq:derivative-kappa} that $\kappa$ is decreasing.
   Since $g$ is superharmonic,  $g'(a+)-g'(a-) \leq 0$ by Proposition \ref{prop:Riesz}; therefore, by \eqref{eq:kappa}, $\kappa(a+)-\kappa(a) \leq 0$ for $a \in (a_1, a_2) \cap F$. Hence, 
    $\kappa$ is decreasing on $(a_1, a_2)$. 
    The converse follows by an analogous argument.
    \end{proof}

\section{Examples}\label{sec:example}
\subsection{One-sided stopping region case}\label{subsection:onesided}
We consider the McKean optimal stopping problem with respect to spectrally negative Lévy process (\cite{mordecki, alili-kyp, christensen2009note}). The problem is given by
\[v_{M}(x) = \sup_{\tau \in \mathcal{T}} \E_x ( e^{-q \tau} g(x)), \ \text{where} \ g_M(x) = (K-e^x)_+.
\]
We exclude the case where $q=0$ but $\lim_{t \goes} X_t  \neq \infty$ with probability $1$ since there is no optimal stopping time in this case. If $q=0$ and $\lim_{t \goes} X_t  = -\infty$ with probability $1$, then  the optimal strategy is to never stop due to $g(-\infty)> g(x)$ for all $x$.  If $q=0$ and $X$ oscillates, then the process is recurrent and there exists no optimal stopping time.
By the exponential change of measure \eqref{eq:expochange}, this problem is equivalent to the following problem:
\( v_{M}(x) = \sup_{\tau \in \mathcal{T}} \E^*_x  [g^*_M(x)], \ \text{where} \  g^*_M (x) = e^{-\Phi(q)x}(K-e^x)_+ \).
Since $g(-\infty) =K$ and $g(x)<K$ for all $x$, $g$ does not attain its maximum.

\begin{figure}[t]
\begin{minipage}{0.48\linewidth}
  \centering
  \includegraphics[width=\linewidth]{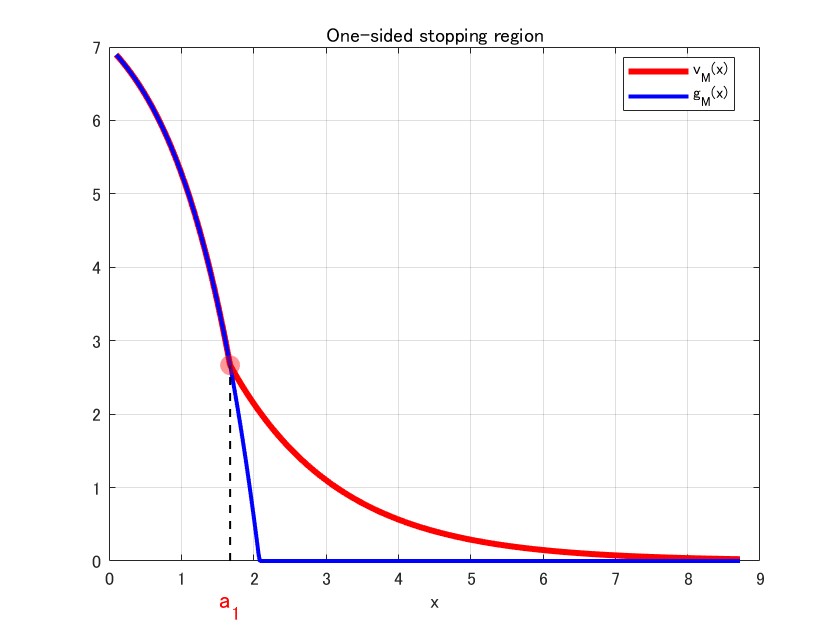}
  \captionsetup{font=scriptsize}
  \caption{One-sided stopping region case}
  \label{fig:onesided}
\end{minipage}\hfill
\begin{minipage}{0.48\linewidth}
  \centering
  \includegraphics[width=\linewidth]{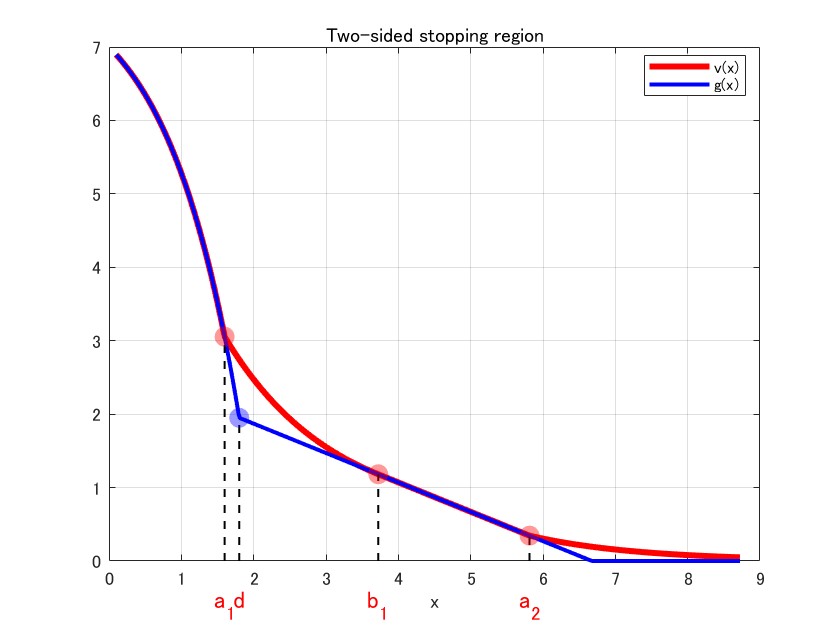}
  \captionsetup{font=scriptsize}
  \caption{Two-sided stopping region case}
  \label{fig:twosided}
\end{minipage}
\end{figure}

Let $\LL^*$ be the generator of $X$ under $\p^*$.
It follows from \eqref{tildeL} and a straightforward calculation that for $x < \log K$,
\begin{align*}
    \LL^* g^*_M (x) &= e^{-\Phi(q)x}(\LL -q) g_M(x)= e^{-\Phi(q)x} \Big( -\psi(1)e^x - q(K-e^x) \Big)<0.
\end{align*}
On the other hand, it follows from \eqref{tildeL} and a straightforward calculation that for $x > \log K$,
\begin{align*}
     \LL^* g^*_M  (x) &= e^{-\Phi(q)x}(\LL -q) g_M(x)  = e^{-\Phi(q)x} \int_{(y<\log K -x)} (K-e^{x+y}) \Pi(\diff y) >0.
\end{align*}
Hence, Corollary \ref{coro:one-sided} applies.
Let $\kappa$ be the $\kappa$-function associated with $g^*_M$.
It follows from a direct calculation that  
\[ \kappa(a) = e^{-\Phi(q)a} \left[K \frac{q}{\Phi(q)}-\frac{q-\psi(1)}{\Phi(q)-1}e^a\right].\]
In particular, if $q=0$, then 
\(\kappa(a) =  K \psi'(0+)-\psi(1)e^a. \)
Hence, by Corollary \ref{coro:one-sided}, we obtain the optimal threshold $a_1$  as follows
\begin{align}
a_1
&=\log \Big(K \frac{q}{\Phi(q)}\frac{\Phi(q)-1}{q-\psi(1)} \Big);\quad q>0,\quad
a_1
=\log \Big(K \psi'(0+) /\psi(1) \Big);\quad q=0.
\label{astar}
\end{align}
Therefore, we obtain that  the value function $v=H_{a_1} g$ and the optimal stopping time $\tau^* = T_{(-\infty, a_1)}$.
The corresponding graph of the solution is shown in Figure~\ref{fig:onesided}.
\footnote{For the numerical analysis, we take $\psi(\theta)=c \theta -\lambda (1-\mu(\mu+\theta)^{-1})$ with $\mu = 1.5$, $\lambda = 1$, and $c = 1.2$.}

\subsection{Two-sided stopping region case}\label{subsection:twosided}
We define the reward function \( g \) as follows:
\[
g(x)
=
(K-e^x)\mathbf{1}_{\{x\le d\}}
+
\max\!\{(K-e^d)-l(x-d),\,0\}\mathbf{1}_{\{x\ge d\}},
\]
where $K>0$, $l>0$, $d < \log K$ are constants.
To simplify the discussion, we assume $q=0$.
If $\lim_{t \goes} X_t = -\infty$ holds with probability $1$ or $X$ oscillates, then there is no optimal stopping time as in Section \ref{subsection:onesided}. Hence, we can assume $\lim_{t \goes} X_t = \infty$ holds with probability $1$ without loss of generality.
Let $v$ be the value function.
Let $a^*$ denote the value of $a_1$ defined in Section \ref{subsection:onesided}.
To exclude trivial cases, we additionally assume that the functions $g$ and $v_{M}$ intersect for $x > d$ and $d \in (a^*, \log K)$, where $v_{M}$ is  defined in Section \ref{subsection:onesided}.
Let $p>d$ be the smallest intersection point of $g$ and $v_{M}$.
Note that as \( d \to \log K \), the problem reduces to the McKean optimal stopping problem in Section \ref{subsection:onesided}.
We can show that $g$ does not attain its maximum  by the same argument in Section \ref{subsection:onesided}

 
Henceforth, we specify $X$ as the process whose Laplace exponent is
\begin{align}
   \psi(\theta)=c \theta -\lambda (1-\mu(\mu+\theta)^{-1}), \label{eq:exopojimps}
\end{align}
where $c$ is the drift rate, $\lambda>0$ is the rate of the arrival rate, and $\mu>0$ is the parameter associated with the exponentially distributed jumps. We assume that $c-\lambda/\mu >0$ so that $\psi'(0)>0$.
Then, the scale function $W$, as given in \cite{Hubalek_Kyprianou_2009}, and $h_a$ are expressed as follows: 
\[W(x) = c^{-1} (1 + \lambda(c \mu - \lambda)^{-1} (1 - e^{-(\mu - c^{-1} \lambda)x} ) ) \conn\] 
\[h_a(x) = K- (\psi(1)/\psi'(0))e^a+\lambda(c \mu - \lambda)^{-1}(\mu + 1)^{-1} e^a e^{-(\mu - c^{-1} \lambda)(x-a)}.\]
We can find $a_1$ and $b_1$ by applying Proposition \ref{prop:verification}.
Define $g^1:= H_{a_1,b_1}g$ and apply Proposition \ref{prop:verification}  to  $g^1$  again.
We can find $a_2 < \infty$ and $b_2=\infty$, and the itteration ends.
Define $E = [a_1, b_1] \cup [a_2, \infty)$, and 
the value function and the optimal stopping time are given by
\[v = H_{\R \setminus E} g \conn \tau^* = T_{\R \setminus E}. \]
The parameters are set as $(\mu, \lambda, c, K , l, d) =(1.5, 1, 1.2, 8, 0.4, 1.8)$.
Under these values, the computed solutions are $a_1 = 1.5986$, $b_1 = 3.7229$, and $a_2 = 5.8136$. 
The corresponding graph of the solution is shown in Figure~\ref{fig:twosided}.

\subsection{A discontinuous reward}\label{subsection:discontinuous}

We consider the following reward function: 
\[
g(x) = (K-e^x)\mathbf{1}_{\{x < c\}}.
\]
where $K>0$ and $c \leq \log K$.
For simplicity, we consider the case where $q=0$.
Note that this reward is not upper semi continuous. However, as seen in Section \ref{sec:sufficient}, we can apply Proposition \ref{prop:verification}
 and Corollary \ref{coro:one-sided}.
 As shown in Section \ref{subsection:onesided}, $\LL g$ is negative on $(-\infty, c)$, and 
it follows from $g=g'=g''=0$ on $(c,\infty)$ that $\LL g$ is positive on $(c,\infty)$.
Since the number of subharmonic components $n=1$, Corollary \ref{coro:one-sided}  applies.
It follows from a direct calculation that
\begin{align*}
\kappa(a)=
    \begin{cases}
         K \psi'(0+)-\psi(1)e^a  &\text{if} \ a \leq c\\
         -\int_0^\infty \int_{(x,\infty)} g(x+a-z) \Pi(- \diff z) \diff x  &\text{if }\ a>c,
    \end{cases}
\end{align*} 
Let $a_1 = \inf \{a: \kappa(a)\leq 0  \}$ and let $N:= \log(K\psi'(0+)/\psi(1))$. Then, we have
\begin{align*}
    a_1 =
    \begin{cases}
        \log (K\psi'(0+)/\psi(1)) &\text{if} \ N \leq c\\
        c &\text{if} \ N > c,
    \end{cases}
\end{align*}
and $\kappa(a_1) = 0$ if $N \leq c$ and $\kappa(a_1) =\kappa(c) = K \psi'(0+)-\psi(1)e^c >0$ if $N > c$.
By Corollary \ref{coro:one-sided}, the value function and the optimal stopping time are given by $v(x) = H_{a_1} g(x)$  and $\tau^* = T_{(-\infty, a_1)}$.
By using \eqref{eq:expectedreward} and the definition of the smooth Gerber-Shiu function, we have 
\begin{align*}
 v(a_1)-v(a_1-)
 &=H_{a_1} g(a_1)-g(a_1-)\\ 
&= h_{a_1}(a_1)-W(0+)\kappa(a_1)-g(a_1-)\\
&=-W(0+)\kappa(a_1).
\end{align*}
Therefore, the continuous fit condition fails if and only if 
$N>c$.
Figure \ref{fig:onesided_NoCFC} presents the value function and the reward function in the case $ N> c$, with the process specified as in \eqref{eq:exopojimps}.

\begin{figure}[t]
\begin{minipage}{0.48\linewidth}
  \centering
  \includegraphics[width=\linewidth]{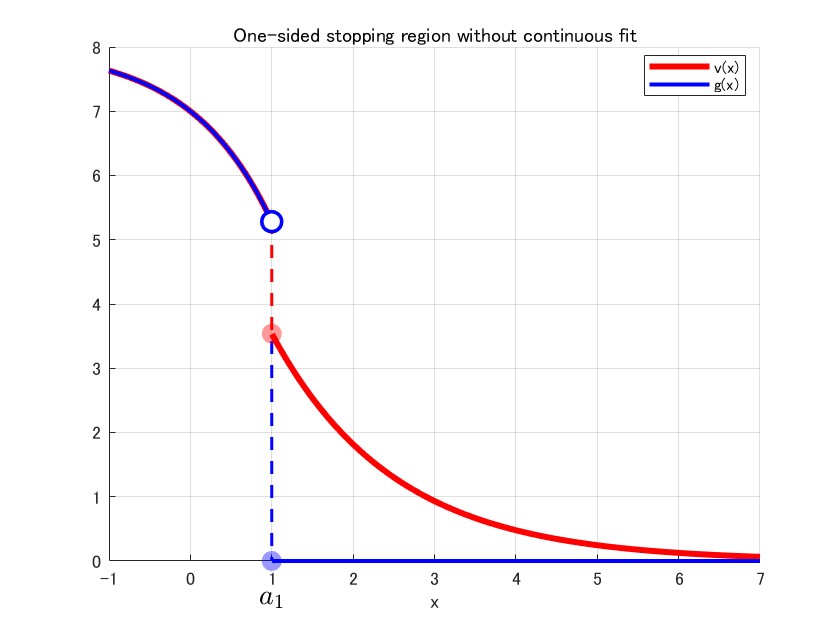}
  \captionsetup{font=scriptsize}
  \caption{One-sided stopping region without continuous fit}
  \label{fig:onesided_NoCFC}
\end{minipage}\hfill
\begin{minipage}{0.48\textwidth}
  \centering
  \includegraphics[width=\linewidth]{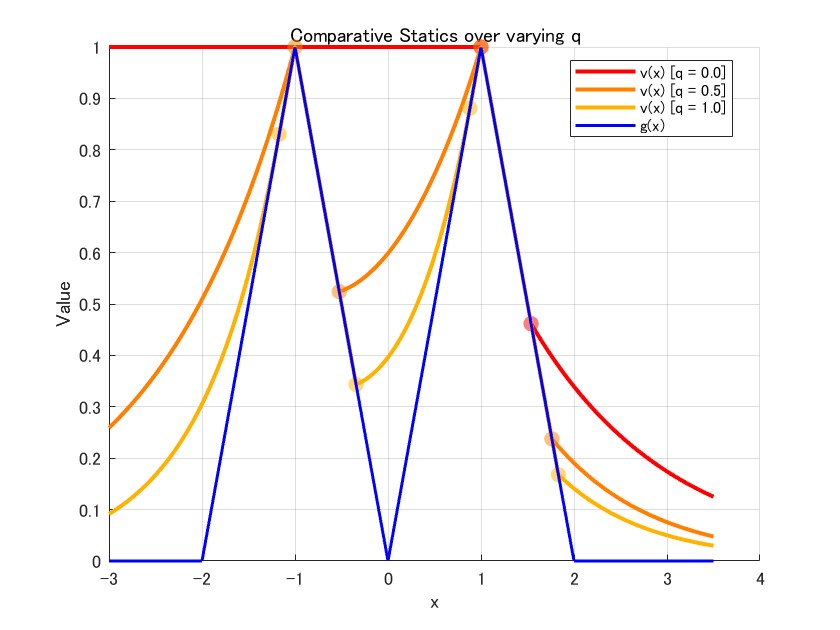}
   \captionsetup{font=scriptsize}
  \caption{Value functions for different discount rates $q=0.0, 0.5$ and $1.0$.}
  \label{fig:Mcurve}
\end{minipage}
\end{figure}

\begin{figure}
    \centering
    \includegraphics[width=0.7\linewidth]{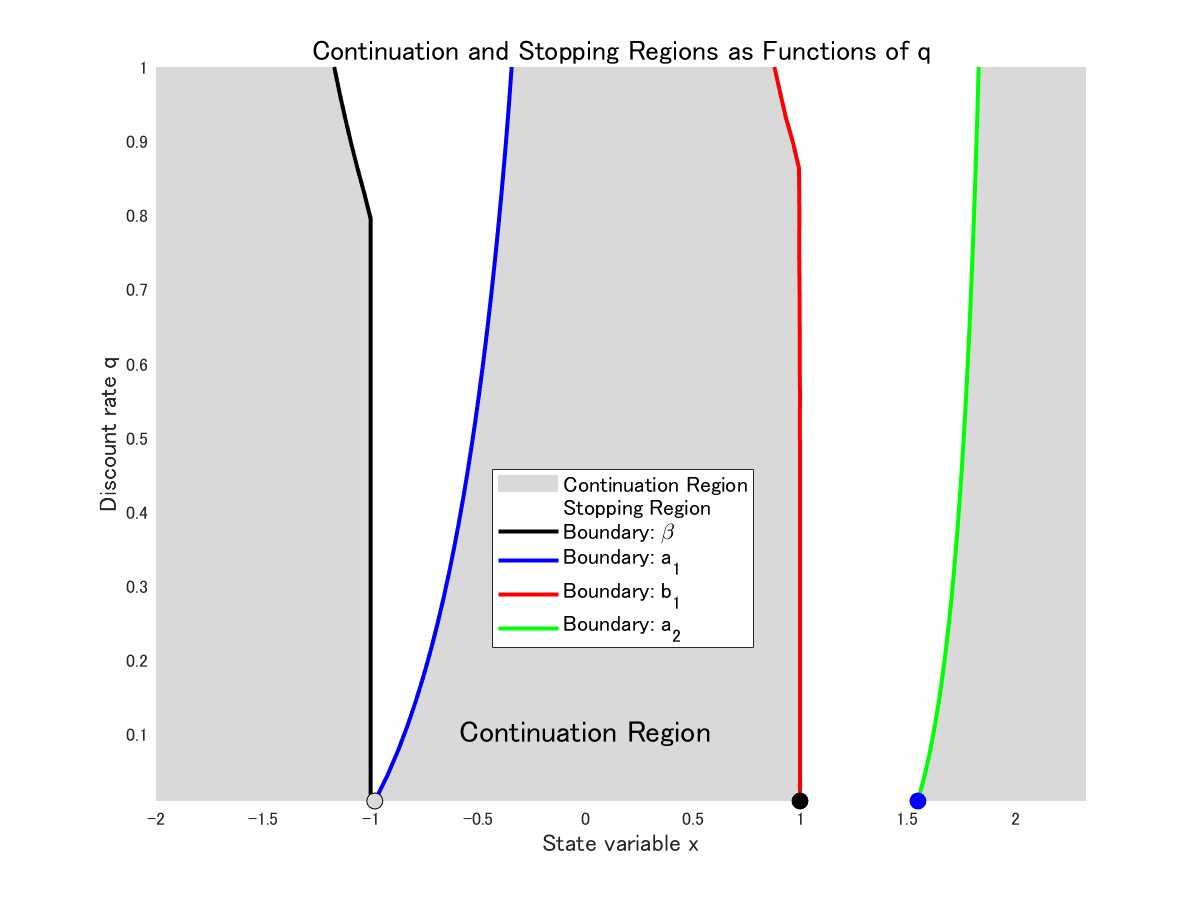}
    \caption{Continuation and stopping regions as functions of $q$}
    \label{fig:Mcurve_conti_q}
\end{figure}

\subsection{A bimodal function case}\label{subsec:bimdal}
We consider the following bimodal reward function: 
\[
g(x)
=
\max(1-|x-1|, \ 0)
+
\max(1-|x+1|, \ 0).
\]
The process is specified as in \eqref{eq:exopojimps} with parameters $(c, \lambda, \mu) = (1.2, 1, 1.5)$.
The structure of the optimal solution depends on the discount rate~$q$.
We present the solution for $q = 0,\ 0.5,$ and $1$, which are numerically solved by applying Corollary  \ref{coro:verification-max}.
Let the value functions and the optimal stopping times be $v_q$ and $\tau^*_q$ for $q=0, 0.5, 1$, respectively.
The solution for $q=0$ is presented in the red curve in Figure \ref{fig:Mcurve}.
The solution exhibits two one-sided continuation regions, 
and it has $m=1$ and $b_1 = \infty$;
\[v_0 = H_{\R \setminus E_0} g \conn \tau^*_0 = T_{\R \setminus E_0} 
\quad \text{with} \quad E_0 = (-\infty, \beta] \cup [a_2, \infty).\]
On the other hand, when $q = 0.5$ or $q=1$, a two-sided continuation region appears in the trough between the two peaks. 
The value function for $q=0.5$ and $q=1$ are presented in the orange and yellow curves in Figures \ref{fig:Mcurve}, respectively .
In the case where $q=0.5$, the solution has $m=2$ with $\beta =-1$, $b_1 = 1$  and $b_2 = \infty$;
\[v_{0.5} = H_{\R \setminus E_{0.5}} g \conn \tau^*_{0.5} = T_{\R \setminus E_{0.5}} 
\quad \text{with} \quad E_{0.5} = (-\infty, -1] \cup [a_1, 1] \cup [a_2, \infty).\]
Note that the smooth fit condition at $b_2=1$ fails.
In the case where  $q=1$, the solution has $m=2$ with $b_1 \in (0, 1)$  and $b_2 = \infty$;
\[v_{1} = H_{\R \setminus E_{1}} g \conn \tau^*_{1} = T_{\R \setminus E_{1}} 
\quad \text{with} \quad E_{1} = (-\infty, \beta] \cup [a_1, b_1] \cup [a_2, \infty).\]

Figure \ref{fig:Mcurve_conti_q} illustrates the comparative statics of the continuation and stopping regions with respect to the discount rate $q$. In this figure, the shaded grey areas represent the continuation region, while the white areas correspond to the stopping regions. The four colored curves denote the boundaries separating these distinct regions: the black line represents boundary $\beta$, the blue line is $a_1$, the red line is $b_1$, and the green line is $a_2$. The primary finding introduced by this graph is that as the discount rate $q$ increases along the vertical axis, the overall continuation region  shrinks, which consequently leads to the expansion of the stopping regions. A notable characteristic of these dynamics is that the black boundary $\beta$ and the red boundary $b_1$ remain pinned exactly at -1 and 1, respectively, for a substantial range of lower discount rates, although at $q = 0$, as confirmed above, $\beta = 1$ and $b_1 = \infty$. They only begin to detach and shift leftward once the discounting becomes sufficiently high. Meanwhile, the blue boundary $a_1$ and the green boundary $a_2$ shift continuously rightward as $q$ increases.

\appendix
\section{Proof of the left continuity of $a \mapsto h_a(x)$}\label{sec:h_a-left-cont}
We give a proof of  the left continuity of $a \mapsto h_a(x)$.

\begin{proposition}\label{prop:h_a-left-cont}
    $a \mapsto h_a(x)$ is left-continuous for each fixed $x$.
\end{proposition}
\begin{proof}
    Fix $x \in \R$ arbitrarily.
    Since $a \mapsto g(a-)$ and $a \mapsto g'(a-)$ are left continuous, we prove the left continuity of $K(\cdot, x)$, where 
    \[K(a,x) := \int_0^{x-a} W(x-a-y) J_a(y) \diff y. \]
    Let $(a_n)$ be an increasing sequence with $a_n \uparrow a$.
    Taking $\delta>0$ sufficiently small, we may assume that $a-a_n<\delta$.
    We have
    \begin{align*}
        \int_0^{x-a_n}& \bigl|W(x-a_n-y) J_{a_n}(y) \bigl|\diff y \\
        &\leq \int_0^{x-a+\delta} W(x-a+\delta) \sup_{a' \in (a-\delta, a), u \in (0, x-a+\delta) } \bigl|J_{a'}(u) \bigl|  \diff y.
    \end{align*}
    Therefore, it suffices to establish \begin{align}
        \sup_{a' \in (a-\delta, a), u \in (0, x-a+\delta) } \bigl|J_{a'}(u) \bigl| < \infty \quad  \text{for some} \quad \delta>0, \label{eq:J-sup}
    \end{align}
    together with the left continuity of $a \mapsto J_a(x)$.
    Indeed, once these are available, the dominated convergence theorem yields the left continuity of $K(\cdot,x)$.
    To prove \eqref{eq:J-sup}, we prepare several lemmas. The left continuity of $a \mapsto J_a(x)$ is obtained as a by-product of the same argument (see the comments after Lemma~\ref{lemma:I-conti}).
\end{proof}

We denote the integral of the absolute value of the jump term in $J_a(x)$ by $I(a,x)$; that is,
\[I(a,x) :=  \int_{(x,\infty)} \big| g(x+a-z)-g(a-)+g'(a-)(z-x) \big| \Pi(-\diff z). \]
Since $g$ is locally bounded, it suffices to show that the following to obtain \eqref{eq:J-sup} ;
\begin{align}
\sup_{a' \in (a-\delta, a), u \in (0, x-a+\delta) } I(a',u) < \infty \quad  \text{for some} \quad \delta>0. \label{eq:I-sup}    
\end{align}

From now on, we fix $a \in \R$ and $x>0$.
Take $\delta>0$ such that $(a-2\delta, a) \cap F = \varnothing$.

\begin{lemma}\label{lemma:I-finite}
    For $(a',u) \in [a-\delta, a] \times [0, x-a+\delta]$,
    \[I(a',u) < \infty.\]
\end{lemma}
\begin{proof}
    We only consider the case $u=0$, as the case $u>0$ is easier since the integral is bounded away from the singularity at $z=0$. Fix $a'$ and take $\eps>0$ such that $(a'-\eps, a') \cap F = \varnothing$.  We have
    \begin{align*}
        I(a',0) 
        \leq& \int_{(0, \eps)} \big|g(a'-z)-g(a'-)+g'(a'-)z \big| \Pi (-\diff z) \\
        &+ \int_{(\eps, \infty)} \big|g(a'-z)-g(a'-)+g'(a'-)z \big| \Pi (-\diff z). 
    \end{align*}
By the Taylor's theorem, the first term of right-hand side is bounded from above by
\[  \frac{1}{2} \sup_{v \in (0,\eps)} |g''(a'-v)|  \int_{(0, \eps)} z^2 \Pi(-\diff z)< \infty.  \]
It follows from Assumptions \ref{assumption-jumpsize} (iv) and (v) that the second term is bounded from above by the following finite integral
\[
\int_{(\eps, \infty)}|g(a'-z)| \Pi(-\diff z) + |g(a'-)|\Pi(-\infty, -\eps)+ |g'(a'-)|\int_{(\eps, \infty)} z \Pi(-\diff z).\]
Theretofore, we conclude $I(a',0) < \infty$.
\end{proof}

\begin{lemma}\label{lemma:I-conti}
     $I$ is jointly continuous on $[a-\delta, a] \times [0, x-a+\delta]$.
\end{lemma}
\begin{proof}
    Let $(a', u) \in [a-\delta, a] \times [0, x-a+\delta]$ and take sequences $(a_n) \subset [a-\delta, a]$ and $x_n \subset [0, x-a+\delta]$ with $a_n \to a'$ and $(x_n) \to u$. We only consider the case $u=0$, as  the case $u>0$ is easier since the integral is bounded away from the singularity at $z=0$.
    We have
    \begin{align*}
        I(a_n, x_n) =& \int_{(x_n, \delta)} \big| g(x_n+a_n-z)-g(a_n-)+g'(a_n-)(z-x_n) \big| \Pi(-\diff z) \\
        &+ \int_{(\delta, \infty)} \big| g(x_n+a_n-z)-g(a_n-)+g'(a_n-)(z-x_n) \big| \Pi(-\diff z),
    \end{align*}
    where we define the first and second term of the right-hand side as $I_1(a_n, x_n)$ and $I_2(a_n, x_n)$, respectively.

    Using Taylor's theorem,
    \begin{align*}
        I_1(a_n, x_n) 
        &\leq \frac{1}{2}\sup_{v \in (0,2\delta)}\big|g''(a-v) \big| \int_{(x_n, \delta)} (z-x_n)^2 \Pi(-\diff z) \\
        &\leq \frac{1}{2}\sup_{v \in (0,2\delta)}\big|g''(a-v) \big| \int_{(0, \delta) }z^2 \Pi(-\diff z)<\infty.
    \end{align*}
    Using Assumptions \ref{assumption-jumpsize} (iv) and (v), $I_2(a_n, x_n)$ is bounded from above by the following finite integral:
    \begin{align*}
        \int_{(\delta,\infty)} \left( \sup_{v \in (x+a-2\delta, x+a+\delta)}|g(v-z)|+ \sup_{v \in (a-\delta, a)}|g(v)|+\sup_{v \in (a-\delta,a)} |g'(v)|z \right)   \Pi(-\diff z).
        \end{align*}
        Therefore we can apply the dominated convergence theorem  and obtain 
        \begin{align*}
           \lim_{n \goes} I(a_n, x_n)
        &= \int_{(x_n, \infty)} \lim_{n \goes} \big| g(x_n+a_n-z)-g(a_n-)+g'(a_n-)(z-x_n) \big| \Pi(-\diff z)\\
        &= I(a',u), 
        \end{align*}
        where the second equality follows from $g$ is $\Pi$-almost every point by Assumption \ref{assumption-jumpsize}  (ii), (iii) and (vi).
        Therefore,   
        $I$ is continuous on $[a-\delta, a] \times [0, x-a+\delta]$.  
\end{proof}

Since $a \mapsto g(a-)$ and $a \mapsto g'(a-)$ are left continuous, Lemma \ref{lemma:I-conti} implies the left continuity of $a \mapsto J_a(x)$ for each fixed $x$.

\begin{lemma}
    \eqref{eq:I-sup} holds. Hence, \eqref{eq:J-sup} also holds.
\end{lemma}
\begin{proof}
    It follows from Lemmas \ref{lemma:I-finite} and \ref{lemma:I-conti} that $I$ is real-valued continuous function on the compact set $[a-\delta, a] \times [0, x-a+\delta]$. Therefore, $I$ has a finite maximum and  \eqref{eq:I-sup} holds.
\end{proof}

\bibliographystyle{abbrv}
{\small \bibliography{references}}

\end{document}

%% file: realline2_modified.tex
\begin{figure}[htbp]
  \centering
  \begin{tikzpicture}[scale=1.2]

    \draw[thick,->] (-5,0) -- (5,0);

    \coordinate (a_0) at (-4,0);
    \coordinate (l1) at (-3,0);
    \coordinate (r1) at (-1.5,0);
    \coordinate (l2) at (-0.4,0);
    \coordinate (r2) at (0.6,0);
    \coordinate (l3) at (2.2,0);
    \coordinate (l0) at (4,0);

    \node[below] at (l1) {$\ell_1$};
    \node[below] at (r1) {$r_1$};
    \node[below] at (l2) {$\ell_2$};
    \node[below] at (r2) {$r_2$};
    \node[below] at (l3) {$\ell_3$};
    \node[below] at (l0) {$\ell_0$};

    \fill (l1) circle (2pt);
    \fill (r1) circle (2pt);
    \fill (l2) circle (2pt);
    \fill (r2) circle (2pt);
    \fill[red] (l3) circle (2pt);
    \fill (l0) circle (2pt);

    \draw[ultra thick, red] (-2.95,0) -- (-1.55,0);
    \node[above, red] at (-2.4, 0.15) {$\LL g > 0$};
    \draw[ultra thick, red] (-0.35,0) -- (0.55,0);
    \node[above, red] at (0, 0.15) {$\LL g > 0$};
    \draw[ultra thick, red] (4.05,0) -- (4.95,0);
    \node[above, red] at (4.5, 0.15) {$\LL g > 0$};
    \node[above, red] at (2.4, 0.15) {$g'(\ell_3+)>g'(\ell_3-)$};


  \end{tikzpicture}
  \captionsetup{font=scriptsize} 
  \caption{The red segments represent the set $D =(\ell_1, r_1)\cup (\ell_2, r_2) \cup \{\ell_3\} \cup (\ell_0, \infty) $, which consists of  four subharmonic components. }
  \label{fig:realline2}
\end{figure}

%% file: potential_July2026.bbl
\begin{thebibliography}{10}

\bibitem{alili-kyp}
L.~Alili and A.~E. Kyprianou.
\newblock Some remarks on first passage of {L}\'{e}vy processes, the {A}merican put and smooth pasting.
\newblock {\em Ann. Appl. Probab.}, 15:2062--2080, 2004.

\bibitem{avram2020w}
F.~Avram, D.~Grahovac, and C.~Vardar-Acar.
\newblock The \textit{W}, \textit{Z} scale functions kit for first passage problems of spectrally negative {{Lévy}} processes, and applications to control problems.
\newblock {\em ESAIM: Probability and Statistics}, 24:454--525, 2020.

\bibitem{avram2015gerber}
F.~Avram, Z.~Palmowski, and M.~R. Pistorius.
\newblock On gerber--shiu functions and optimal dividend distribution for a {{Lévy}} risk process in the presence of a penalty function.
\newblock {\em Ann. Appl. Probab.}, 25 (4):1868 -- 1935, 2015.

\bibitem{Bertoin_1996}
J.~Bertoin.
\newblock {\em {L}\'evy processes}, volume 121 of {\em Cambridge Tracts in Mathematics}.
\newblock Cambridge University Press, Cambridge, 1996.

\bibitem{Bertoin_1997}
J.~Bertoin.
\newblock Exponential decay and ergodicity of completely asymmetric {L}\'evy processes in a finite interval.
\newblock {\em Ann. Appl. Probab.}, 7:156--169, 1997.

\bibitem{biffis2010note}
E.~Biffis and A.~E. Kyprianou.
\newblock A note on scale functions and the time value of ruin for {{Lévy}} insurance risk processes.
\newblock {\em Insurance: Mathematics and Economics}, 46(1):85--91, 2010.

\bibitem{Chan_2009}
T.~Chan, A.~Kyprianou, and M.~Savov.
\newblock Smoothness of scale functions for spectrally negative {L}\'evy processes.
\newblock {\em Probab. Theory Relat. Fields}, 150(3--4):691--708, 2011.

\bibitem{christensen2019optimal}
S.~Christensen, F.~Crocce, E.~Mordecki, and P.~Salminen.
\newblock On optimal stopping of multidimensional diffusions.
\newblock {\em Stochastic Processes and their Applications}, 129(7):2561--2581, 2019.

\bibitem{christensen2009note}
S.~Christensen and A.~Irle.
\newblock A note on pasting conditions for the american perpetual optimal stopping problem.
\newblock {\em Statistics \& probability letters}, 79(3):349--353, 2009.

\bibitem{christensen2018multidimensional}
S.~Christensen and P.~Salminen.
\newblock Multidimensional investment problem.
\newblock {\em Mathematics and Financial Economics}, 12:75--95, 2018.

\bibitem{christensen2013optimal}
S.~Christensen, P.~Salminen, and B.~Q. Ta.
\newblock Optimal stopping of strong markov processes.
\newblock {\em Stochastic Processes and their Applications}, 123(3):1138--1159, 2013.

\bibitem{deligiannidis2009optimal}
G.~Deligiannidis, H.~Le, and S.~Utev.
\newblock Optimal stopping for processes with independent increments, and applications.
\newblock {\em Journal of applied probability}, 46(4):1130--1145, 2009.

\bibitem{egami-yamazaki2013}
M.~Egami and K.~Yamazaki.
\newblock Precautionary measures for credit risk management in jump models.
\newblock {\em Stochastics}, 85(1):111--143, 2013.

\bibitem{Egami-Yamazaki-2011}
M.~Egami and K.~Yamazaki.
\newblock On the continuous and smooth fit principle for optimal stopping problems in spectrally negative {{Lévy}} models.
\newblock {\em Adv. in Appl. Probab.}, 46:139--167, 2014.

\bibitem{Egami_Yamazaki_2010_2}
M.~Egami and K.~Yamazaki.
\newblock Phase-type fitting of scale functions for spectrally negative {L}\'evy processes.
\newblock {\em Journal of Computational and Applied Mathematics}, 264:1--22, 2014.

\bibitem{hsiau2014logconcave}
S.-R. Hsiau, Y.-S. Lin, and Y.-C. Yao.
\newblock Logconcave reward functions and optimal stopping rules of threshold form.
\newblock 2014.

\bibitem{Hubalek_Kyprianou_2009}
F.~Hubalek and A.~Kyprianou.
\newblock Old and new examples of scale functions for spectrally negative \lev processes.
\newblock {\em Sixth Seminar on Stochastic Analysis, Random Fields and Applications, eds R. Dalang, M. Dozzi, F. Russo. Progress in Probability, Birkhäuser}, 2011.

\bibitem{karatzas}
I.~Karatzas and S.~E. Shreve.
\newblock {\em Brownian Motion and Stochastic Calculus}.
\newblock Springer Science\(+\)Business Media, New York, 2nd edition, 1998.

\bibitem{kyprianou2013gerber}
A.~E. Kyprianou.
\newblock {\em Gerber--Shiu risk theory}.
\newblock Springer Science \& Business Media, 2013.

\bibitem{Kyprianou_2014}
A.~E. Kyprianou.
\newblock {\em Fluctuations of {L}\'evy {P}rocesses with {A}pplications}.
\newblock Universitext. Springer-Verlag, Berlin, 2nd edition, 2014.

\bibitem{kyprianou2005novikov}
A.~E. Kyprianou and B.~A. Surya.
\newblock On the {N}ovikov-{S}hiryaev optimal stopping problems in continuous time.
\newblock 2005.

\bibitem{Kyprianou_Surya_2007}
A.~E. Kyprianou and B.~A. Surya.
\newblock Principles of smooth and continuous fit in the determination of endogenous bankruptcy levels.
\newblock {\em Finance Stoch.}, 11:131--152, 2007.

\bibitem{lin2019one}
Y.-S. Lin and Y.-C. Yao.
\newblock One-sided solutions for optimal stopping problems with logconcave reward functions.
\newblock {\em Advances in Applied Probability}, 51(1):87--115, 2019.

\bibitem{mordecki}
E.~Mordecki.
\newblock Optimal stopping and perpetual options for {L}\'{e}vy processes.
\newblock {\em Finance and Stochastics}, 6:473--493, 2002.

\bibitem{mordecki2016optimal}
E.~Mordecki and Y.~Mishura.
\newblock Optimal stopping for {{Lévy}} processes with one-sided solutions.
\newblock {\em SIAM Journal on Control and Optimization}, 54(5):2553--2567, 2016.

\bibitem{mordecki2021two}
E.~Mordecki and F.~Oli{\'u}~Eguren.
\newblock Two-sided optimal stopping for {{Lévy}} processes.
\newblock 2021.

\bibitem{mordecki2007optimal}
E.~Mordecki and P.~Salminen.
\newblock Optimal stopping of hunt and {{Lévy}} processes.
\newblock {\em Stochastics An International Journal of Probability and Stochastic Processes}, 79(3-4):233--251, 2007.

\bibitem{novikov2007solution}
A.~Novikov and A.~Shiryaev.
\newblock On a solution of the optimal stopping problem for processes with independent increments.
\newblock {\em Stochastics An International Journal of Probability and Stochastic Processes}, 79(3-4):393--406, 2007.

\bibitem{protter2012stochastic}
P.~E. Protter.
\newblock Stochastic integration and differential equations.
\newblock Springer, 2012.

\bibitem{rw-2}
L.~Rogers and D.~Williams.
\newblock {\em Diffusions, Markov processes and martingales, Vol.2: 2nd Edition}.
\newblock Cambridge University Press, Cambridge, 1994.

\bibitem{salminen1985}
P.~Salminen.
\newblock Optimal stopping of one-dimensional diffusions.
\newblock {\em Math. Nachr.}, 124(1):85--101, 1985.

\bibitem{Shiryaev_2008}
A.~N. Shiryaev.
\newblock {\em Optimal stopping rules}, volume~8 of {\em Stochastic Modelling and Applied Probability}.
\newblock Springer-Verlag, Berlin, 2008.

\bibitem{surya2007approach}
B.~A. Surya.
\newblock An approach for solving perpetual optimal stopping problems driven by {{Lévy}} processes.
\newblock {\em Stochastics An International Journal of Probability and Stochastic Processes}, 79(3-4):337--361, 2007.

\bibitem{Surya_2008}
B.~A. Surya.
\newblock Evaluating scale functions of spectrally negative {L}\'evy processes.
\newblock {\em J. Appl. Probab.}, 45:135--149, 2008.

\end{thebibliography}
